\documentclass[12pt]{article}
\usepackage{fullpage}
\usepackage{amsfonts, amsmath, amsthm,amssymb}

\newtheorem{theorem}{Theorem}[subsection]
\newtheorem{lemma}[theorem]{Lemma}
\newtheorem{prop}[theorem]{Proposition}
\newtheorem{cor}[theorem]{Corollary}
\newtheorem{question}[theorem]{Question}
\theoremstyle{definition}
\newtheorem{defn}[theorem]{Definition}
\newtheorem{example}[theorem]{Example}
\newtheorem{remark}[theorem]{Remark}

\newtheorem{hypo}[theorem]{Hypothesis}
\newtheorem{notation}[theorem]{Notation}

\numberwithin{equation}{theorem}

\def\calE{\mathcal{E}}
\def\calR{\mathcal{R}}
\def\calS{\mathcal{S}}
\def\tcalR{\tilde{\mathcal{R}}}

\def\CC{\mathbb{C}}
\def\QQ{\mathbb{Q}}
\def\RR{\mathbb{R}}
\def\ZZ{\mathbb{Z}}
\def\dual{\vee}
\def\gothm{\mathfrak{m}}
\def\gotho{\mathfrak{o}}
\def\be{\mathbf{e}}
\def\bv{\mathbf{v}}
\def\bw{\mathbf{w}}
\def\bx{\mathbf{x}}

\DeclareMathOperator{\alg}{alg}
\DeclareMathOperator{\an}{an}
\DeclareMathOperator{\bd}{bd}
\DeclareMathOperator{\coker}{coker}
\DeclareMathOperator{\con}{con}
\DeclareMathOperator{\Ext}{Ext}
\DeclareMathOperator{\Frac}{Frac}
\DeclareMathOperator{\GL}{GL}
\DeclareMathOperator{\Hom}{Hom}
\DeclareMathOperator{\inte}{int}
\DeclareMathOperator{\rank}{rank}
\DeclareMathOperator{\tr}{tr}

\begin{document}

\title{Slope filtrations for relative Frobenius}
\author{Kiran S. Kedlaya}
\date{September 6, 2007}

\maketitle

\begin{abstract}
The slope filtration theorem gives a partial analogue
of the eigenspace decomposition of a linear transformation,
for 
a Frobenius-semilinear endomorphism of a finite free module over the 
Robba ring (the ring 
of germs of rigid analytic functions on an unspecified open annulus of
outer radius 1) over a discretely valued field.
In this paper, we give a third-generation proof of this theorem, which 
both introduces some new simplifications (particularly the use of faithfully
flat descent, to recover the theorem from a classification theorem of
Dieudonn\'e-Manin type) and extends the result to allow an arbitrary
action on coefficients (previously the action on coefficients had to itself
be a lift of an absolute Frobenius).
This extension is relevant to a study of $(\phi, \Gamma)$-modules associated
to families of $p$-adic Galois representations, as initiated
by Berger and Colmez.
\end{abstract}

\section*{Introduction}

This paper describes a third-generation proof of the
slope filtration theorem for Frobenius modules over the Robba ring
(Theorem~\ref{T:slope filt} herein).
This proof is more expedient than what one finds
in our original paper \cite{kedlaya-local} or its sequel
\cite{kedlaya-slope}. 
In addition, we 
generalize the slope filtration theorem by allowing for ring endomorphisms
which do not act as Frobenius lifts on scalars, only on the series
variable. This is intended as a prelude to a theory of Frobenius
modules in families; we will not develop such a theory here, but 
see the next section for reasons one might want to do so, from the
realm of $p$-adic Hodge theory.
(Note that \cite{kedlaya-slope} itself generalizes \cite{kedlaya-local} in
a different direction, replacing the power series rings by somewhat more 
general
objects; we do not treat that generalization here.)

For
an alternate perspective on this theorem and some related results in
$p$-adic differential equations and $p$-adic Hodge theory,
we also recommend Colmez's Bourbaki notes \cite{colmez-bourbaki}.

\subsection{Context}

The slope filtration theorem
\cite[Theorem~6.10]{kedlaya-local} (also exposed in \cite{kedlaya-slope})
gives a partial classification of
Frobenius-semilinear transformation on finite free
modules over the Robba ring (a certain ring of univariate formal Laurent
series with $p$-adic coefficients).
It is loosely analogous to the
eigenspace decomposition of a linear transformation in ordinary
linear algebra; it is also closely related to Manin's classification
of rational Dieudonn\'e modules.

The slope filtration theorem was originally introduced in the context of 
Berthelot's rigid
cohomology, a $p$-adic Weil cohomology
for varieties in characteristic $p$.
There, one obtains a analogue of the $\ell$-adic local monodromy theorem, 
originally conjectured by Crew \cite{crew};
this analogue can be used to establish
various structural results such as finiteness of cohomology
\cite{kedlaya-finite} and purity in the sense of Deligne \cite{kedlaya-weilii}. 

The effect of the slope 
filtration theorem on $p$-adic Hodge theory has perhaps
been even more acute: it enables one
to study $p$-adic Galois representations via their
associated $(\phi,\Gamma)$-modules over the Robba ring. This point of view
has been put forth chiefly by Berger with striking consequences:
he has proved Fontaine's conjecture that de Rham representations
are potentially semistable \cite{berger-cst}, and given
an alternate proof of the Colmez-Fontaine theorem
on admissibility of filtered $(\phi,N)$-modules \cite{berger-weak}.
(A useful variant of the latter argument has been given by Kisin
\cite{kisin}.)
More recently Colmez \cite{colmez-tri} used this viewpoint to define a class
of \emph{trianguline representations} of a $p$-adic Galois group;
these play an important role in the $p$-adic
local Langlands correspondence for $\GL_2(\QQ_p)$ \cite{colmez-local}.
The trianguline are also important in the theory of $p$-adic modular forms, as
most local Galois representations attached to overconvergent $p$-adic
modular forms (namely, those of noncritical slope)
 are trianguline.
The $p$-adic local Langlands correspondence in turn
has touched off a flurry of activity, which this
introduction is not the right place to summarize; we merely note
the resolution of Serre's conjecture by
Khare-Wintenberger \cite{kw, kw2},
and progress on the
Fontaine-Mazur conjecture by Kisin \cite{kisin-fm} 
and Emerton (in preparation).

In both rigid cohomology and $p$-adic Hodge theory,
one is led to study Frobenius modules in families,
i.e., over the Robba ring with coefficients not in a $p$-adic field but
in, say, an affinoid algebra. 
In either situation, the first step to studying
Frobenius modules in families is to pass from a family to a generic point,
which on rings amounts to replacing an integral affinoid algebra with a
complete field containing it.
In the rigid cohomology version of this
argument, the resulting field is itself acted on by Frobenius, so the
slope filtration theorem as
presented in \cite{kedlaya-local, kedlaya-slope} 
is immediately applicable; indeed, the key technique in
\cite{kedlaya-finite} is to extend the application of the local monodromy
theorem on the generic point to a large enough subspace of the base space. 
However, in the
$p$-adic Hodge theory version, one might like to allow ``Frobenius'' to
act in some fashion on the base of the family other than simply a lift of
the $p$-power map; in fact, one natural situation is where the base is not
moved at all.

One goal of this paper, and in fact the principal reason for its existence,
 is to generalize the slope filtration theorem to
modules over the Robba ring with an action of a ``relative Frobenius'', which
may do whatever one wishes to coefficients as long as it acts like a
Frobenius lift on the series parameter. We hope this will lead to some
study of $p$-adic Hodge theory in families; some of the corresponding
analysis in equal characteristics has been initiated by Hartl \cite{hartl},
using an equal-characteristic analogue of the slope filtration theorem
based on the work of Hartl and Pink \cite{hartl-pink}. 
In mixed characteristics,
Hartl \cite{hartl2} has set up part of a corresponding theory,
which addresses a conjecture of Rapoport and Zink \cite{rapoport-zink}
from their work on period spaces for $p$-divisible groups; results are
presently quite fragmentary, but a good theory of $(\phi, \Gamma)$-modules
in families may help. Another potential application would be to
analysis of the local geometry of the Coleman-Mazur eigencurve
\cite{coleman-mazur}, which parametrizes the Galois representations
attached to certain $p$-adic modular forms, or of higher-dimensional
``eigenvarieties'' associated to automorphic representations on 
groups besides $\GL_2$. An initial step in this direction has already
been taken by Bella\"\i che-Chenevier \cite{bellaiche-chenevier},
who study deformations of trianguline representations; however, this
involves only a zero-dimensional base, so they can already apply the
usual slope filtration theory after
a restriction of scalars. For other questions, e.g., properness, one would
want to consider positive-dimensional bases like a punctured disc.
In this direction,
Berger and Colmez have introduced
a theory of \'etale $(\phi, \Gamma)$-modules associated
to $p$-adic Galois representations in families
\cite{berger-colmez},
which relativizes some of the results of Cherbonnier-Colmez 
\cite{cherbonnier-colmez} and 
Berger \cite{berger-weak} for a single $p$-adic Galois representation. 

\subsection{About the results}

For the sake of introduction, we give here a very brief description of 
what the original slope filtration theorem says, how the main result of
this paper extends it, and what novelties in the argument are
introduced in this paper. Start with a complete
discretely valued field $K$ of mixed characteristics $(0,p)$.
Let $\calR$ be the ring of formal Laurent series $\sum_{n \in \ZZ} c_n u^n$
convergent on some annulus with outer radius 1 (but whose inner radius may
depend on which series is being considered). Let $\phi_K: K \to K$ be
an endomorphism lifting the absolute $q$-power Frobenius 
on the residue field
of $K$, for some power $q$ of $p$,
and define a map $\phi: \calR \to \calR$ by the formula
$\phi(\sum c_n u^n) = \sum \phi_K(c_n) \phi(u)^n$,
where $\phi(u) - u^q$ has all coefficients of norm less than 1.
Let $M$ be a finite free $\calR$-module equipped with a
$\phi$-semilinear map $F: M \to M$ which takes any basis of $M$ to another basis
of $M$ (it is enough to check for a single basis). Then 
\cite[Theorem~6.10]{kedlaya-slope} asserts that
$M$ admits an exhaustive 
filtration whose successive quotients are each pure of
some slope (i.e., some power of $F$ times some scalar acts on some
basis via an invertible matrix over the subring of $\calR$ of series with
\emph{integral} coefficients), 
and the slopes increase as you go up the filtration; moreover,
those requirements uniquely characterize the filtration.

As noted earlier, 
the slope filtration should be thought of as analogous to what one might get
from a linear transformation over $K$ by grouping eigenspaces, interpreting
the slope of an eigenspace as the valuation of its eigenvalue. One can
in fact deduce 
an analogous such result for semilinear transformations over $K$,
which also follows from the Dieudonn\'e-Manin classification theorem.
One might then expect that the slope filtration can be generalized so
as to allow \emph{any} isometric
action on $K$, not just a Frobenius lift;
that is what is established in this paper
(Theorem~\ref{T:slope filt}).

As promised earlier in this introduction, one happy side effect of this
generalization is the introduction of
some technical simplifications.
We give a development of the theory of slopes which does not
depend on already having established the Dieudonn\'e-Manin-style 
classification; this follows up on a suggestion made in
\cite{kedlaya-slope}. We give a much simplified version of the descent 
argument that deduces the filtration theorem from the DM classification,
based on the idea of replacing the Galois descent used previously
with faithfully flat descent; this avoids the use of comparison between
generic and special Newton polygons, and of some intricate approximation
arguments. (In particular, there is no longer any need to deal with 
finite extensions of the Robba ring, which allows for some notational
and expository simplifications.)
That substitution creates some flexibility in what we may take
as the ``extended Robba ring'' for the DM classification; here we use
a ring made from generalized power series, some of whose properties are 
a bit more transparent than for the corresponding ``big rings'' in
\cite{kedlaya-local} and \cite{kedlaya-slope}.

\subsection{Structure of the paper}

The structure of this paper is a bit unusual, as we have attempted
to make the paper more friendly to the novice reader by fronting some of the
key assertions and pushing back more technical aspects. (This
assertion applies both to the paper as a whole, and to
Sections~\ref{sec:dm} and~\ref{sec:descend} individually.)
The consequence is that
the logical structure is a bit loopy: results are stated, and sometimes
used, before having been proved. However, we hope that it is not too hard
to see that there are indeed no vicious circles in the reasoning.

In Section~\ref{sec:filt thm}, we introduce the Robba ring, the category
of $\phi$-modules, the notions of degree and slope, the subcategories
of pure $\phi$-modules of various slopes, and the statement of the
filtration theorem.

In Section~\ref{sec:dm}, we introduce an extended Robba ring (whose elements
are modeled on Hahn-Mal'cev-Neumann generalized power series rather than
ordinary power series), state a classification theorem for $\phi$-modules
over the extended Robba ring, then perform the calculations required to prove
this theorem.

In Section~\ref{sec:descend}, we deduce the slope filtration theorem from
the classification theorem over the extended Robba ring. The key tool here
is an invocation of faithfully flat descent for modules.

\subsection*{Acknowledgments}
Thanks to Laurent Berger for the original suggestion to consider
relative Frobenius and for subsequent
discussions, to Lucia di Vizio for providing the reference to
Praagman's work, and to Peter Schneider for additional comments. 
The author was supported by NSF grant DMS-0400727,
NSF CAREER grant DMS-0545904, and a Sloan Research Fellowship.

\section{Statement of the filtration theorem}
\label{sec:filt thm}

\subsection{The Robba ring}

\begin{defn} \label{D:initial}
Let $K$ be a field complete for a discrete valuation,
with residue field $k$; let $\gotho_K$ denote the
valuation subring of $K$ and let $\gothm_K$ denote the maximal ideal of 
$\gotho_K$. 
(We need not make any restriction on the characteristics of $K, k$.)
Write $|\cdot|$ for some fixed norm corresponding
to the valuation (the normalization does not matter).
For $r>0$, let $\calR^r$ be the ring of rigid analytic functions on the
annulus $e^{-r} \leq |t| < 1$ (these are just Laurent series in
the variable $t$ convergent on this
region), and let $\calR$ be the union of the $\calR^r$. The ring $\calR$
is called the \emph{Robba ring} over $K$.
It follows from the work of Lazard 
\cite{lazard}
that $\calR$ is a \emph{B\'ezout domain},
that is, an integral domain in which every finitely generated ideal is principal.
\end{defn}

\begin{remark}
Any B\'ezout domain $R$ enjoys a number of nice properties generalizing
properties of principal ideal domains,
including the following.
Some of these are actually properties of \emph{Pr\"ufer domains}, in which
every finitely generated ideal is projective; these generalize Dedekind
domains to the non-noetherian setting.
\begin{itemize}
\item
Any finite locally free $R$-module is free
\cite[Proposition~2.5]{kedlaya-local}.
\item
Any torsion-free $R$-module is flat; this holds for any Pr\"ufer domain
 \cite[VII Proposition~4.2]{cartan-eilenberg}.
\item
Any finitely presented projective $R$-module is free
\cite[Proposition~4.8]{crew}.
\item
If $M$ is a finite free $R$-module and $N$ is a submodule of $M$ which is 
saturated,
i.e., $N = M \cap (N \otimes_R \Frac R)$, then $N$ and $M/N$ are both free
\cite[Proposition~4.8]{crew}, \cite[Lemma~2.4]{kedlaya-local}.
\end{itemize}
\end{remark}

\begin{defn}
Let $\calR^{\inte}$ be the subring of $\calR$ consisting of series with 
coefficients in $\gotho_K$; this ring is a discrete valuation ring with
residue field $k ((t ))$, which is not complete but is 
henselian \cite[Lemma~3.9]{kedlaya-local}. 
Let $\calR^{\bd}$ be the subring of $\calR$ consisting
of series with bounded coefficients; it is the fraction field 
of $\calR^{\inte}$.
\end{defn}

\begin{remark} \label{R:bounded by 1}
Note that for $x \in \calR$, one has
$x \in \calR^{\inte}$ if and only if there exists an integer $n$
such that the function $t^n x$ is bounded by 1 on some annulus
$e^{-r} \leq |t| < 1$. 
\end{remark}

\begin{remark} \label{R:same units}
Lazard's work \cite{lazard}
includes a theory of Newton polygons for elements of $\calR$,
using which one can read off numerous structural properties. One key
example is that
the units in $\calR$ are precisely the nonzero elements of $\calR^{\bd}$
\cite[Corollary~3.23]{kedlaya-local}.
\end{remark}

\begin{remark}
One can also define the Robba ring even if the valuation on $K$ is not
discrete, but its properties are very different. For instance,
$\calR^{\bd}$ is no longer the fraction field of $\calR^{\inte}$.
This makes even the formulation of a slope theory over such $K$,
let alone any proofs, somewhat more delicate than the approach we take here.
\end{remark}

\subsection{Frobenius lifts on the Robba ring}

\begin{defn}
Fix an integer $q>1$. (To see why we forbid $q=1$, see 
Remark~\ref{R:no identity}.)
A \emph{relative ($q$-power)
Frobenius lift} on the Robba ring is a homomorphism
$\phi: \calR \to \calR$ of the form $\sum_i c_i t^i \mapsto \sum_i \phi_K(c_i)
u^{i}$, where $\phi_K$ is an isometric field 
endomorphism of $K$ and $u \in \calR^{\inte}$ is such that $u - t^q$ is in the
maximal ideal of $\calR^{\inte}$. If $k$ has characteristic $p>0$ and
$q$ is a power of $p$,
we define an \emph{absolute ($q$-power)
Frobenius lift} as a
relative Frobenius lift in which $\phi_K$ is itself a $q$-power
Frobenius lift. 
\end{defn}

\begin{remark}
The treatments in \cite{kedlaya-local, kedlaya-slope} only allow
absolute Frobenius lifts, and the approaches do not carry over easily to the
general case because of the use of Galois descent at some key moments.
See the introduction for discussion of why one needs the relative case.
\end{remark}

\begin{defn}
For $r > 0$, let $|\cdot|_r$ denote the supremum norm on the 
circle $|t| = e^{-r}$, as applied to elements of $\calR^r$; one
easily verifies that
\[
\left| \sum_{i \in \ZZ} c_i t^i \right|_r = \sup_i \{|c_i| e^{-ri}\}.
\]
We extend the definition to vectors by taking the maximum over entries.
\end{defn}

\begin{remark} \label{R:entire}
Note that for $f$ analytic on the entire open unit disc (i.e., represented
by an ordinary power series rather than a Laurent series), 
we have $|f|_r \leq |f|_s$
whenever $0 < s \leq r$; in other words, the supremum of $f$ over the entire
disc $|t| \leq e^{-s}$ occurs on the circle $|t| = e^{-s}$.
In fancier language, the circle $|t| = e^{-s}$ is the \emph{Shilov boundary}
of the disc $|t| \leq e^{-s}$, as in \cite[Corollary~2.4.5]{berkovich}.
\end{remark}

\begin{remark} \label{R:annuli}
Let $\phi$ be a relative Frobenius lift; then for some $r_0 > 0$, we have
$|\phi(t)/t^q - 1|_{r_0/q} < 1$. It follows that for $r \in (0, r_0)$ and
$f \in \calR^r$, $\phi(f) \in \calR^{r/q}$ and
$|f|_r = |\phi(f)|_{r/q}$.
In geometric terms, $\phi$ induces a surjective map
from the annulus $e^{-r/q} < |t| < 1$ to the annulus $e^{-r} < |t| < 1$.
(Compare \cite[Lemma~3.7]{kedlaya-local}.)
\end{remark}

The following is both a typical example of how to make calculations on
Robba rings and a crucial ingredient in what follows.
\begin{prop} \label{P:h1 both ways}
Let $\phi$ be a relative Frobenius lift, and
let $A$ be an $n \times n$ matrix over $\calR^{\inte}$. Then the map
$\bv \mapsto \bv - A \phi(\bv)$ on column vectors induces a bijection
on $(\calR/\calR^{\bd})^n$.
\end{prop}
\begin{proof}
The problem is unaffected if we replace $\bv, A$ by $t^m \bv, (t^m/\phi(t^m))
A$, so by Remark~\ref{R:bounded by 1}, we may
reduce to the case where the entries of $A$ are bounded by 1 
on some annulus with outer
radius 1.
Choose $r_0$ as in Remark~\ref{R:annuli}.
To check injectivity, we must argue that if 
$\bw = \bv - A \phi(\bv)$ is bounded,
then so is $\bv$. 
Choose
$r \in (0,r_0)$ such that $A, \bw, \phi(\bv)$ have entries which are defined
on the annulus $e^{-r} \leq |t| < 1$, and the entries of $A$ are bounded 
by 1 there. 
Choose $c>0$ such that $|\bw|_s \leq c$ for $0<s \leq r$, and such that
$|\phi(\bv)|_s \leq c$ for 
$r/q \leq s \leq r$. (The latter is possible because every
analytic function on a closed annulus is bounded.)
Then $|\bv|_s = |\bw + A \phi(\bv)|_s \leq c$ for $r/q \leq s \leq r$,
so $|\phi(\bv)|_s \leq c$ for $r/q^2 \leq s \leq r/q$.
Repeating the argument, we see that $|\bv|_s \leq c$ for $0 < s \leq r$,
proving the claim.
(Compare \cite[Lemma~3.3.3]{kedlaya-slope}.)

To check surjectivity, take $\bw \in \calR^n$. 
Choose
$r \in (0,r_0)$ such that $A, \bw$ have entries which are defined
on the annulus $e^{-r} \leq |t| < 1$, and the entries of $A$ are bounded by 1 there. 
Define the sequence $\{\bw_l\}_{l=0}^\infty$ as follows. Start with
$\bw_0 = \bw$. Given $\bw_l$, write $\bw_l = \sum_{i \in \ZZ} \bw_{l,i} t^i$,
put $\bw_l^+ = \sum_{i>0} \bw_{l,i} t^i$ and $\bw_l^- = \bw_l - \bw_l^+$,
and put $\bw_{l+1} = A \phi(\bw_l^+)$.
Since the entries of $t^{-1} \bw_l^+$ are analytic on the entire open unit 
disc,
by Remark~\ref{R:entire} we have 
\[
|\bw_l^+|_{r} \leq e^{-r+r/q} |\bw_l^+|_{r/q} \leq e^{-r+r/q} |\bw_l|_{r/q};
\]
consequently, $|\bw_{l+1}|_{r/q} \leq e^{-r+r/q} |\bw_l|_{r/q}$.
Thus the sequence $\bw_l^+$
converges to zero under $|\cdot|_{r/q}$, and hence also
under $|\cdot|_s$ for $s \geq r/q$ by Remark~\ref{R:entire}.
On the other hand, for $0 < s \leq r/q$, 
applying Remark~\ref{R:entire} after substituting $t \mapsto t^{-1}$ gives
\[
|\bw_l^-|_s \leq |\bw_l^-|_{r/q} \leq |\bw_l|_{r/q}.
\]

Now set $\bv = \sum_{l=0}^\infty \bw_l^+$; then $\bv$ has entries analytic on
the closed disc of radius $e^{-r/q}$, and $\bw - \bv + A \phi(\bv) = \sum_{l=0}^\infty
\bw_l^-$ is bounded on $e^{-r/q} \leq |t| < 1$. Since $\phi(\bv)$ is
analytic on the closed disc of radius $e^{-r/q^2}$, we can write
$\bv = \bw + A \phi(\bv) - \sum_{l=0}^\infty \bw_l^-$ 
and thus extend $\bv$ across the annulus
$e^{-r/q} \leq |t| \leq e^{-r/q^2}$; by induction, $\bv$ extends to
the entire open unit disc.
This proves the desired surjectivity.
\end{proof}

One can also prove the following, as in \cite[Lemma~5.4.1]{kedlaya-slope}.
\begin{prop} \label{P:cherbonnier1}
Let $\calE$ denote the $\gothm_K$-adic completion of
$\calR^{\bd}$.
Let $\phi$ be a relative Frobenius lift on $\calR$, and 
let $A$ be an $n \times n$ matrix over $\calR^{\inte}$. If $\bv
\in \calE^n$ is a column vector such that
$A \bv = \phi (\bv)$, then $\bv \in (\calR^{\bd})^n$.
\end{prop}
\begin{proof}
This will follow later from Proposition~\ref{P:cherbonnier2};
we will not use it in the interim.
\end{proof}
\begin{remark} \label{R:cherbonnier1}
In the case where $A$ is invertible, 
Proposition~\ref{P:cherbonnier1}
was proved independently by Cherbonnier (unpublished,
but see \cite[Th\'eor\`eme~III.1.1]{cherbonnier-colmez})
and Tsuzuki \cite[Proposition~4.1.1]{tsuzuki-etale}. Tsuzuki's
underlying argument 
can be used even when $A$ is not invertible; see 
\cite[Proposition~2.2.2]{tsuzuki-etale}.
\end{remark}

\begin{remark}
It should be possible to carry everything in this paper over to the
case where one only assumes $\phi(t) = \sum_i c_i t^i$ such that 
$c_q \in \gotho_K^*$ and $c_i \in \gothm_K$ for $i < q$.
(For instance, in the theory of $(\phi, \Gamma)$-modules, the composition
of the usual $\phi$ with any nontrivial $\gamma \in \Gamma$ would have this
property.) The proof of Proposition~\ref{P:h1 both ways} extends to this
setting, but the embedding of $\calR$ into the extended Robba ring
$\tcalR$ of Section~\ref{sec:dm} must be modified, as accordingly
must the projection construction of Section~\ref{sec:descend}.
\end{remark}

\subsection{$\phi$-modules}

\begin{defn}
Define a \emph{$\phi$-(ring/field)}
 to be a ring/field $R$
equipped with an endomorphism $\phi$; we say $R$ is \emph{inversive}
if $\phi$ is bijective.
Define a \emph{(strict) $\phi$-module} over a $\phi$-ring $R$
to be a finite free $R$-module $M$
equipped with an isomorphism $\phi^* M \to M$, which we also think of as a
semilinear $\phi$-action on $M$; the semilinearity means that for $r \in R$
and $m \in M$, $\phi(rm) = \phi(r) \phi(m)$. 
Note that the category of $\phi$-modules admits tensor products, 
symmetric and exterior
powers, and duals.
\end{defn}

\begin{remark}
The definition of $\phi$-module used here is somewhat more restrictive than
one sees in other contexts, hence the optional modifier ``strict''. 
For instance, in some cases one allows
modules which are projective but not free, or worse.
In other cases, one allows the $\phi$-action to take kernel and cokernel in
some $\phi$-stable Serre category of $R$-modules; we will do this ourselves
shortly.
\end{remark}

\begin{remark}
It will be convenient for us to describe $\phi$-modules in terms of bases
and matrices. If $M$ is a $\phi$-module and $\be_1, \dots, \be_n$ is a
basis of $M$, we can completely describe the $\phi$-action on $M$ by specifying
the invertible $n \times n$ matrix $A$ which satisfies $\phi(\be_j) = \sum_i
A_{ij} \be_i$. Note that the semilinearity skews conjugation:
if $\be'_1, \dots, \be'_n$ is another basis and
the change of basis matrix $U$ is defined by $\be'_j = \sum_i U_{ij} \be_i$,
then the $\phi$-action on the new basis is via the matrix
$U^{-1} A \phi(U)$.
\end{remark}

It is also useful to think of
$\phi$-modules as modules for a twisted polynomial ring.
\begin{defn}
For $R$ a $\phi$-ring,
define the \emph{twisted polynomial ring} $R\{T\}$ to be the set of
finite formal sums $\sum_{i=0}^\infty a_i T^i$ with $a_i \in R$,
equipped with the noncommutative
ring structure in which $T a = \phi(a) T$ for $a \in R$.
If $R$ is a field, then all left ideals 
of $R\{T\}$ are principal, by the division algorithm \cite[Theorem~6]{ore}.
If $R$ is inversive, 
one may similarly define a \emph{twisted Laurent polynomial ring}
$R\{T^{\pm}\}$.
\end{defn}

\begin{remark}
In general, a $\phi$-module over $R$ can be interpreted as a
left $R\{T\}$-module which is finite free over $R$, but one must remember
the condition that $\phi$ carries some basis to another basis. 
On the other hand, if $R$ is inversive,
then the data of a $\phi$-module over $R$ is equivalent to the data
of a left $R\{T^{\pm}\}$-module which is finite free over $R$.
If $R$ is an inversive $\phi$-field, then
irreducible $\phi$-modules over $R$ all have the form
$R\{T^{\pm}\}/R\{T^{\pm}\}P$ for some irreducible twisted polynomial $P$.
\end{remark}

When talking about pure slopes, it will be helpful to switch from
working with $\phi$ to working with a power of $\phi$; the following 
definition facilitates this switch.
\begin{defn} \label{D:pull push}
View $\phi$-modules as left modules for the twisted polynomial ring
$R\{T\}$. For $a$ a positive integer, define the
\emph{$a$-pushforward functor} $[a]_*$ from $\phi$-modules to $\phi^a$-modules
to be the restriction along the inclusion $R\{T^a\} \to R\{T\}$.
Define the \emph{$a$-pullback functor} $[a]^*$ from $\phi^a$-modules
to $\phi$-modules to be the extension of scalars functor
\[
M \mapsto R\{T\} \otimes_{R\{T^a\}} M.
\]
The following are easily verified (as in
\cite[\S 3.2]{kedlaya-slope}):
\begin{itemize}
\item
The functors $[a]^*$ and $[a]_*$ form an adjoint pair.
\item
The functors $[a]_*$ and $[a]^*$ are exact and commute with duals; consequently,
$[a]_*$ and $[a]^*$ also form an adjoint pair (i.e., in the other order).
\item
The functor $[a]_*$ commutes with tensor products over $R$ (but
$[a]^*$ does not).
\item
If $M$ is a $\phi$-module and $N$ is a $\phi^a$-module, then
$M \otimes [a]^* N \cong [a]^* ([a]_* M \otimes N)$.
\item
If $M$ is a $\phi$-module, then $\rank([a]_* M) = \rank(M)$.
\item
If $N$ is a $\phi^a$-module, then $\rank([a]^* N) = a \rank(N)$.
\item
If $N$ is a $\phi^a$-module, then
$[a]_* [a]^* N \cong N \oplus \phi^*(N) \oplus 
\cdots \oplus (\phi^{a-1})^*(N)$.
\end{itemize}
\end{defn}

\begin{defn} \label{D:h0 h1}
For $M$ a $\phi$-module, put
\[
H^0(M) = \ker(\phi-1: M \to M), \qquad
H^1(M) = \coker(\phi-1: M \to M).
\]
One easily checks that in the category of $\phi$-modules over $R$,
\[
\Hom(M, N) \cong H^0(M^\dual \otimes N), \qquad
\Ext(M, N) \cong H^1(M^\dual \otimes N).
\]
Moreover, for $N$ a $\phi^a$-module, there are natural bijections
\[
H^i(N) \cong H^i([a]^* N) \qquad (i=0,1).
\]
\end{defn}

\begin{remark}
Beware that although the pullback/pushforward terminology was inspired
by a related construction in \cite{hartl-pink}, the two do not agree in
that context.
\end{remark}

\subsection{Degrees, slopes, and stability}

For the rest of this section, we will put ourselves in the following 
situation. Note that Hypothesis~\ref{hypo:robba} 
has a weak form and a strong form; we will
assume only the weak form unless otherwise specified.
(Thanks to Peter Schneider for suggesting this dichotomy.)
\begin{hypo} \label{hypo:robba}
Let $R^{\inte} \subseteq R^{\bd} \subseteq R$ be inclusions of B\'ezout domains
such that $R^* \subset R^{\bd}$.
Let $\phi$ be an endomorphism of $R$ acting also on $R^{\bd}$ and $R^{\inte}$.
Let $w: R^{\bd} \to \ZZ \cup \{+\infty\}$ be a $\phi$-equivariant
valuation such that $w(R^*) = \ZZ$ and 
$R^{\inte} = \{r \in R^{\bd}: w(r) \geq 0\}$.
Suppose in addition that for any $n \times n$ matrix $A$ over $R^{\inte}$,
the map $\bv \mapsto \bv - A \phi(\bv)$ on column vectors induces an injection
(weak form) or bijection (strong form)
on $(R/R^{\bd})^n$. Note that the analogous hypothesis for $\phi^a$ also holds,
since one can identify the kernel and cokernel of $\bv \mapsto \bv - A
\phi^a(\bv)$ on $(R/R^{\bd})^n$ with the kernel and cokernel of
\[
(\bv_0, \bv_1, \dots, \bv_{a-1}) \mapsto (\bv_0 - A \phi(\bv_{a-1}),
\bv_1 - \phi(\bv_0), \dots, \bv_{a-1} - \phi(\bv_{a-2}))
\]
on $(R/R^{\bd})^{na}$. (Compare the last remark in Definition~\ref{D:h0 h1}.)
\end{hypo}

\begin{example} \label{exa:Robba}
For our purposes, the principal example of strong 
Hypothesis~\ref{hypo:robba} is
as follows. We take $R, R^{\bd}, R^{\inte} = \calR, \calR^{\bd}, \calR^{\inte}$
to be the Robba ring and variants over $K$; note that $\calR^{\bd} =
\calR^* \cup \{0\}$.
We take $\phi$ to be a relative Frobenius lift, and $w$ to be
the valuation on $\calR^{\bd}$
for which $\calR^{\inte}$ is the valuation subring.
The last condition in strong
Hypothesis~\ref{hypo:robba} holds by virtue of 
Proposition~\ref{P:h1 both ways}.
We will construct a variation of this example, the extended Robba ring $\tcalR$,
in Section~\ref{sec:dm}; using the axiomatic approach avoids some repetition.
\end{example}

\begin{example} \label{ex:hartl-pink}
Besides the Robba ring, additional examples of strong Hypothesis~\ref{hypo:robba}
are also possible.
Here is one from the work of Hartl and Pink \cite{hartl-pink}:
take $\CC$ to be the
completed algebraic closure of a local field of equal characteristic $p$, 
$R$ to be the Laurent series over $\CC$
convergent on the punctured open unit disc,
$R^{\bd}$ to be the series which are meromorphic at zero,
$\phi$ to be the map $\sum c_i t^i \mapsto \sum c_i^q t^i$ for $q$
a power of $p$, and
$w$ to be the order of vanishing at 0. 
See Remark~\ref{R:hartl-pink} 
and Question~\ref{Q:q-analogue} for further discussion around this example.
\end{example}

\begin{defn}
For $M$ a $\phi$-module over $R$ of rank $n$, the top
exterior power $\wedge^n M$ has rank 1 over $R$; let $\bv$ be a generator,
and write $\phi(\bv) = r\bv$ for some $r \in R^*$. Define 
the \emph{degree} of $M$ by setting $\deg(M) = w(r)$; note that this does not
depend on the choice of the generator 
by virtue of the $\phi$-equivariance of $w$.
If $M$ is nonzero, define the \emph{slope} of $M$ by setting 
$\mu(M) = \deg(M)/\rank(M)$.
\end{defn}

\begin{remark} \label{R:pull push}
Keeping in mind that degree is analogous to the valuation of the determinant (of
a linear transformation on a finite dimensional vector space over a valued field),
the following formal properties are easily verified
(as in \cite[\S 3.4]{kedlaya-slope}).
\begin{itemize}
\item
If $0 \to M_1 \to M \to M_2 \to 0$ is exact, then
$\deg(M) = \deg(M_1) + \deg(M_2)$; hence $\mu(M)$ is a weighted average
of $\mu(M_1)$ and $\mu(M_2)$.
\item
We have $\mu(M_1 \otimes M_2) = \mu(M_1) + \mu(M_2)$.
\item
We have $\mu(\wedge^i M) = i \mu(M)$.
\item
We have $\deg(M^\dual) = -\deg(M)$ and $\mu(M^\dual) = - \mu(M)$.
\item
If $M$ is a $\phi$-module, then $\mu([a]_* M) = a \mu(M)$.
\item
If $N$ is a $\phi^a$-module, then $\mu([a]^* N) = a^{-1} \mu(N)$.
\end{itemize}
\end{remark}

By analogy with the theory of vector bundles, we make the following definition.
\begin{defn}
We say a $\phi$-module $M$ is \emph{(module-)semistable}
if for any nontrivial $\phi$-submodule $N$, we have $\mu(N) \geq \mu(M)$.
We say $M$ is \emph{(module-)stable}
if for any proper nontrivial $\phi$-submodule $N$, we have $\mu(N) > \mu(M)$.
Note that both properties are preserved under \emph{twisting}
(tensoring with a rank 1 module).
\end{defn}

\begin{remark}
In \cite{kedlaya-slope}, the terms ``stable'' and ``semistable'' were
used without the ``module'' modifier; here we will usually retain the
modifier in statements and drop it in proofs.
The modifier is meant to emphasize the difference between this notion of
semistability and the concept of a
``semistable $(\phi, \Gamma)$-module'' in the sense of $p$-adic Hodge theory,
meaning one which appears to come from a semistable Galois representation.
In the end, over the Robba ring the term ``module-semistable'' will be shown to
be synonymous with ``pure'', so the terminological overload will cease to
be a problem.
\end{remark}

\begin{remark}
Those familiar with stability of vector bundles (or with
\cite{hartl-pink}) will notice that our definitions
differ from the usual convention
 by an overall minus sign. The sign convention here (which is also the
one used in \cite{kedlaya-local, kedlaya-slope}) seems to be
more consistent with usage in the theory of crystalline cohomology.
\end{remark}

\begin{prop} \label{P:rank 1 stable}
Any $\phi$-module of rank $1$ is module-stable.
\end{prop}
\begin{proof}
This is a consequence of the assumptions built into weak
Hypothesis~\ref{hypo:robba}.
Namely, by twisting, it suffices to show that the trivial $\phi$-module 
$M \cong R$ is
stable. If $N$ is a $\phi$-submodule of $M$, 
we may write $N = Rx$ for some $x \in M$
such that $\lambda = \phi(x)/x \in R^*$, and by definition
$\mu(N) = w(\lambda)$. If $\mu(N) \leq 0$, then
$x - \lambda^{-1} \phi(x) = 0$ implies $x \in R^{\bd}$ by
weak Hypothesis~\ref{hypo:robba}; hence $N = M$ and $\mu(N) = w(\phi(x)) - w(x) 
= 0$. In other words, $\mu(N) > 0$ unless
$N = M$, as desired.
\end{proof}
\begin{cor} \label{C:same rank}
If $N \subseteq M$ is an inclusion of $\phi$-modules of the same rank, then
$\mu(N) \geq \mu(M)$, with equality if and only if $N = M$.
\end{cor}
\begin{proof}
Put $n = \rank M$
and apply Proposition~\ref{P:rank 1 stable} to the inclusion $\wedge^n N 
\subseteq \wedge^n M$.
\end{proof}

\begin{lemma} \label{L:least slope}
Let $M$ be a $\phi$-module over $R$. Then the slopes of nonzero 
$\phi$-submodules of $M$
are bounded below.
\end{lemma}
\begin{proof}
We proceed by induction on $\rank(M)$. By Corollary~\ref{C:same rank}, the slopes
of $\phi$-submodules of $M$ of full rank are bounded below by $\mu(M)$.
If $M$ has no nontrivial $\phi$-submodules of lower rank, 
then there is nothing more to check.
Otherwise, let $N$ be a saturated $\phi$-submodule of lower rank; then by
hypothesis, the slopes of nonzero $\phi$-submodules of both $N$ and $M/N$ are
bounded below. If now $P$ is any nonzero $\phi$-submodule of $M$, then the sequence
\[
0 \to N \cap P \to P \to P/(N \cap P) \to 0
\]
is exact. If both factors
are nonzero, we have $\mu(N \cap P) \geq \mu(N)$ 
and $\mu(P/(N \cap P)) \geq \mu(M/N)$,
and $\mu(P)$ is a weighted average of $\mu(N \cap P)$ and $\mu(P/(N \cap P))$, so
it is bounded below. If one factor vanishes, then $\mu(P)$ simply equals the slope
of the other factor, so the same conclusion holds.
\end{proof}

\begin{lemma} \label{L:least slope2}
Let $M$ be a nonzero $\phi$-module over $R$. Then there is a largest
$\phi$-submodule of $M$ of least slope, which is module-semistable.
\end{lemma}
\begin{proof}
The fact that there is a least slope $s$ holds by Lemma~\ref{L:least slope}
and the fact that the denominators of slopes are bounded above by the rank
of $M$; clearly any $\phi$-submodule of slope $s$ must be semistable.
If $N_1$ and $N_2$ are two such submodules, 
then the kernel of the surjection $N_1 \oplus N_2
\to N_1 + N_2$ must have slope at least $s$, so $\mu(N_1 + N_2) \leq s$. On the other hand,
$\mu(N_1 + N_2) \geq s$ because $N_1 + N_2 \subseteq M$, so $\mu(N_1 + N_2) = s$.
Hence the $\phi$-submodules of $M$ of slope $s$ are closed under sum, yielding
the existence of a largest such submodule.
\end{proof}
\begin{cor} \label{C:semi push}
Let $M$ be a $\phi$-module over $R$. Then for any positive integer $a$,
$M$ is module-semistable if and only if $[a]_* M$ is module-semistable.
\end{cor}
\begin{proof}
If $[a]_* M$ is semistable, evidently $M$ is too. Conversely, if
$[a]_* M$ is not semistable, then its largest $\phi^a$-submodule of least slope
is a $\phi^a$-submodule $M_1$ of lower rank. By the uniqueness in
Lemma~\ref{L:least slope2}, $M_1$ must in fact be
preserved by $\phi$, so $M$ is not semistable either.
\end{proof}

\begin{defn}
Let $M$ be a $\phi$-module over $R$. 
A \emph{module-semistable filtration} of $M$ is a 
filtration $0 = M_0 \subset M_1 \subset \cdots \subset M_l = M$ by saturated
$\phi$-submodules such that each quotient $M_i/M_{i-1}$ is module-semistable. A 
\emph{Harder-Narasimhan (HN) filtration} is a module-semistable filtration in which
\[
\mu(M_1/M_0) < \cdots < \mu(M_l/M_{l-1}).
\]
\end{defn}

\begin{prop}
Every $\phi$-module over $R$ admits a unique HN filtration, whose first step
is the submodule defined in Lemma~\ref{L:least slope2}.
\end{prop}
\begin{proof}
This is a formal consequence of Lemma~\ref{L:least slope2}; see
\cite[Proposition~4.2.5]{kedlaya-slope}.
\end{proof}

\begin{defn}
Define the \emph{slope multiset} of a module-semistable filtration of a
$\phi$-module of $M$ as the multiset in which each slope of a successive
quotient occurs with multiplicity equal to the rank of that quotient.
These assemble into the lower boundary of a convex region in
the $xy$-plane as follows: start at $(0,0)$, then take each slope $s$
in increasing order and append to the polygon a segment with slope $s$
and width equal to the multiplicity of $s$. The result is called the
\emph{slope polygon} of the filtration; for the HN filtration,
we call the result the \emph{HN polygon}.
\end{defn}
\begin{prop} \label{P:on or above}
The HN polygon lies on or above the slope
polygon of any module-semistable filtration, with the same endpoint.
\end{prop}
\begin{proof}
This is a formal consequence of the definition of an HN filtration: see
\cite[Proposition~3.5.4]{kedlaya-slope}.
\end{proof}

\begin{prop}
Let $M_1, M_2$ be $\phi$-modules over $R$ such that each slope of the
HN polygon of $M_1$ is less than each slope of the HN polygon of $M_2$.
Then $\Hom(M_1,M_2) = 0$.
\end{prop}
\begin{proof}
Choose $f \in \Hom(M_1,M_2)$.
Let $N_1$ be the first step in the HN filtration of $M_1$; then
either $f(N_1) = 0$, or $\mu(f(N_1)) \leq \mu(N_1)$. The latter is impossible
because $\mu(f(N_1))$ is no less than the least slope of $M_2$, whereas
$\mu(N_1)$ is no greater than the greatest slope of $M_1$. Hence
$f$ factors through $M_1/N_1$; repeating,
we obtain $f = 0$.
\end{proof}

\subsection{\'Etale $\phi$-modules}

\begin{defn}
A $\phi$-module $M$ over $R$ or $R^{\bd}$
is said to be \emph{\'etale} (or \emph{unit-root})
if it can be obtained by base extension from
a (strict) $\phi$-module over $R^{\inte}$; that is, $M$ must admit an 
$R^{\inte}$-lattice $N$ such that
$\phi$ induces an isomorphism $\phi^* N \to N$. We call such an $N$ an \emph{\'etale
lattice} of $M$. Note that $N$ is not in general
unique; for instance, it may be rescaled. Note also that the dual of
an \'etale $\phi$-module is again \'etale.
\end{defn}

\begin{remark}
The term ``unit-root'' is standard in applications to crystalline cohomology, where
it refers to the process of extracting the unit roots (roots of valuation 0)
of a $p$-adic polynomial. By contrast,
the term ``\'etale'' is standard in applications to
$p$-adic Hodge theory.
\end{remark}

One of the basic results about \'etale $\phi$-modules is that in a certain
sense, they do not lose
information when base-changed from $R^{\bd}$ to $R$. This can be deduced from
a slightly more general result, which we already used once (to justify that
the Robba ring satisfies Hypothesis~\ref{hypo:robba}) and will use again shortly
(in the proof of Theorem~\ref{T:pure semi}).

\begin{defn}
Define an \emph{isogeny $\phi$-module} over $R^{\inte}$
to be a finite free $R^{\inte}$-module $M$ equipped with an injection
$\phi^* M \to M$ whose cokernel is killed by some power of a uniformizer of
$R^{\inte}$.
Such an object becomes a strict $\phi$-module upon tensoring with
$R^{\bd}$ or $R$.
\end{defn}

\begin{prop} \label{P:h0}
Let $M$ be an isogeny $\phi$-module over $R^{\inte}$.
Then the natural maps $H^i(M \otimes R^{\bd}) \to H^i(M \otimes R)$ for
$i=0$ (under weak Hypothesis~\ref{hypo:robba})
or $i=0,1$ (under strong Hypothesis~\ref{hypo:robba}) are bijective.
\end{prop}
\begin{proof}
This is an immediate consequence of the 
final clause of Hypothesis~\ref{hypo:robba}.
\end{proof}

\begin{prop} \label{P:etale equiv}
The base change functor from \'etale $\phi$-modules over $R^{\bd}$
to \'etale $\phi$-modules over $R$ is an equivalence of categories.
\end{prop}
\begin{proof}
The essential surjectivity holds by definition, so we need only check
full faithfulness. That is, for any \'etale $\phi$-modules $M_1, M_2$
over $R^{\bd}$, we must check that the natural map 
\[
H^0(M_1^\dual \otimes M_2) \to H^0(M_1^\dual \otimes M_2 \otimes R)
\]
is a bijection; this follows from Proposition~\ref{P:h0}.
\end{proof}

\begin{prop} \label{P:etale extend}
Let $0 \to M_1 \to M \to M_2 \to 0$ be a short exact sequence of $\phi$-modules
over $R$. If any two of $M_1, M_2, M$ are \'etale
 (except
possibly $M_1,M_2$ in the case of weak Hypothesis~\ref{hypo:robba}), 
then so is the third.
\end{prop}
\begin{proof}
First, suppose that $M$ and $M_2$ are \'etale. 
By Proposition~\ref{P:etale equiv},
the $\phi$-modules $M, M_2$ and the morphism $M \to M_2$ all descend to 
$R^{\bd}$. By Lemma~\ref{L:etale lattice} below,
we can then produce an \'etale lattice in $M_1$ by taking the
kernel of the map from an \'etale lattice of $M$ to $M_2$.

Next, suppose that $M$ and $M_1$ are \'etale. We then dualize to obtain a second
exact sequence in which $M^\dual$ and $M_1^\dual$ are \'etale. By the previous paragraph,
$M_2^\dual$ is then \'etale, as then is $M_2$.

Finally, suppose that $M_1$ and $M_2$ are \'etale
and that strong Hypothesis~\ref{hypo:robba} holds.
By applying Proposition~\ref{P:h0},
$M_1$, $M_2$, and the exact sequence $0 \to M_1 \to M \to M_2 \to 0$ all descend
to $R^{\bd}$; by rescaling appropriately, we can descend the sequence to
$R^{\inte}$. We can then produce an \'etale lattice in $M$ by lifting
an \'etale lattice from $M_2$, then adding an
\'etale lattice from $M_1$.
\end{proof}

\begin{lemma} \label{L:etale lattice}
Let $M$ be an \'etale $\phi$-module over $R^{\bd}$. Then any finitely generated
$\phi$-stable $R^{\inte}$-submodule of $M$ is a $\phi$-module over $R^{\inte}$.
\end{lemma}
\begin{proof}
Let $M_0$ be an \'etale lattice of $M$, and let $N$ be a finitely generated
$\phi$-stable $R^{\inte}$-submodule of $M$; by rescaling, we may assume
$N \subseteq M_0$. Then $N$ is already an isogeny $\phi$-module,
and it suffices to check that $\deg(N) = 0$; we may do this after
replacing $M$ by $\wedge^{\rank(N)} M$, i.e., we may assume $\rank(N) = 1$.
Let $\be_1, \dots,\be_n$ be a basis of $M_0$, let
$\bv = \sum_{i=1}^n c_i \be_i$ be a generator of $N$,
and put $\phi(\bv) = \sum_{i=1}^n d_i \be_i$. Then
$\deg(N) = \min_i \{w(d_i)\} - \min_i\{w(c_i)\}$, but this
difference is zero because $M_0$ is an \'etale lattice.
\end{proof}

We can also show that \'etale $\phi$-modules are module-semistable, but it will be convenient to do that more generally for pure $\phi$-modules in the
next subsection.

\subsection{Pure $\phi$-modules}

\begin{defn} \label{D:pure}
Let $M$ be a $\phi$-module over $R^{\bd}$ or $R$
of slope $s = c/d$, where $c,d$ are coprime integers
with $d>0$. We say $M$ is \emph{pure} (or \emph{isoclinic}, or sometimes
\emph{isocline}) of slope $s$
if for some $\phi$-module $N$ of rank $1$ and degree $-c$,
$([d]_* M)\otimes N$ is \'etale (the same then holds for any such $N$). 
It will follow from Lemma~\ref{L:pure push} below that
it is equivalent to impose this condition for any one pair $c,d \in \ZZ$ with
$s = c/d$ and $d>0$.
Note that:
\begin{itemize}
\item
any $\phi$-module of rank 1 is pure;
\item
a $\phi$-module is pure of slope 0 if and only if it is \'etale;
\item
the dual of a pure $\phi$-module of slope $s$ is itself pure of slope $-s$.
\end{itemize}
\end{defn}

\begin{remark}
This definition is not that of 
\cite[Definition~6.3.1]{kedlaya-slope}, but it is equivalent to
it by \cite[Proposition~6.3.5]{kedlaya-slope}. It has the advantage that it
can be stated without reference to any sort of Dieudonn\'e-Manin
classification; the downside is that one must expend a bit of effort to check
some natural-looking properties, as we do below.
\end{remark}

\begin{lemma} \label{L:pure push}
Let $M$ be a $\phi$-module over $R^{\bd}$ or $R$, and let $a$
be a positive integer. Then $M$ is pure of some slope $s$ if and only if 
$[a]_* M$ is pure of slope $as$.
\end{lemma}
\begin{proof}
We first check the case where $s=0$. If $M$ is \'etale, then clearly
$[a]_* M$ is too. Conversely, if $[a]_* M$ is \'etale, 
then $\phi$ induces isomorphisms $(\phi^{i+1})^* [a]_* M \to (\phi^{i})^* [a]_* M$
over $R$; by Proposition~\ref{P:etale equiv}, these isomorphisms
descend to $R^{\bd}$. That is, we may reduce to working over $R^{\bd}$.
In this case, let $N_0$ be an \'etale lattice of $[a]_* M$.
Let $N$ be the $R^{\inte}$-span of 
$N_0, \phi(N_0), \dots, \phi^{a-1}(N_0)$; then $N$ is an \'etale
lattice of $M$.
Hence $M$ is \'etale.

In the general case, write $s = c/d$ in lowest terms, and put $b = 
\gcd(a,d)$; then in lowest terms, $as = (ac/b)/(d/b)$. 
Let $N$ be a $\phi^d$-module of rank 1 and degree $-c$;
then $[a/b]_* N$ has rank 1 and degree $-ac/b$.
The following are equivalent:
\begin{itemize}
\item
$M$ is pure of slope $s$;
\item
$([d]_* M) \otimes N$ is \'etale (definition);
\item
$[a/b]_* (([d]_* M) \otimes N) \cong
([ad/b]_* M) \otimes ([a/b]_* N) \cong
([d/b]_* ([a]_* M)) \otimes ([a/b]_* N)$ is \'etale (by above);
\item
$[a]_* M$ is pure of slope $as$ (definition).
\end{itemize}
This yields the claim.
%
\end{proof}
\begin{cor} \label{C:tensor pure}
If $M_1, M_2$ are pure $\phi$-modules of slopes $s_1, s_2$, then
$M_1 \otimes M_2$ is pure of slope $s_1+s_2$.
\end{cor}
\begin{proof}
  By Lemma~\ref{L:pure push}, we may reduce to the case where $s_1, s_2 \in \ZZ$.
By twisting, we may then reduce to the case where $s_1 = s_2 = 0$. In this case the
result follows from the fact that $\phi$-modules over
$R^{\inte}$ admit tensor products.
\end{proof}

We can thus generalize Propositions~\ref{P:etale equiv} and~\ref{P:etale extend}
as follows.
\begin{theorem} \label{T:pure equiv}
For any rational number $s$,
the base change functor from pure $\phi$-modules of slope $s$ over $R^{\bd}$
to pure $\phi$-modules of slope $s$ over $R$ is an equivalence of categories.
\end{theorem}
\begin{proof}
If $M_1, M_2$ are pure of slope $s$, then $M_1^\dual \otimes M_2$ is \'etale.
Hence the proof of Proposition~\ref{P:etale equiv} goes through unchanged.
\end{proof}
\begin{theorem}\label{T:pure ext}
Let $0 \to M_1 \to M \to M_2 \to 0$ be a short exact sequence of 
$\phi$-modules over $R$. If any two of $M_1, M_2, M$ are pure of 
slope $s$ (except
possibly $M_1,M_2$ in the case of weak Hypothesis~\ref{hypo:robba}), 
then so is the third.
\end{theorem}
\begin{proof}
By Lemma~\ref{L:pure push}, we may apply $[a]_*$ to reduce to the case where
$s \in \ZZ$; by twisting, we may force $s=0$. The result now follows from
Proposition~\ref{P:etale extend}.
\end{proof}

\begin{remark}
In a short exact sequence $0 \to M_1 \to M \to M_2 \to 0$ over $R$,
the fact that $M$ is pure of slope $s$ does not by itself
imply the same for $M_1$
and $M_2$, unless the sequence splits (see Corollary~\ref{C:split exact}). 
For example, if $M$ is pure of rank 2 and slope $0$, it can happen
that $M_1$ is pure of rank 1 and slope 1, while $M_2$ is pure of rank 1
and slope $-1$. 
This sort of example arises naturally from $p$-adic Hodge theory,
as in the theory of trianguline representations introduced by Colmez
\cite{colmez-tri}.
\end{remark}

\begin{lemma} \label{L:h0 pos slope}
Let $M$ be a pure $\phi$-module over $R$ of positive slope. Then
$H^0(M) = 0$.
\end{lemma}
\begin{proof}
By replacing $M$ with $[a]_* M$ for $a = \rank(M)$, we can reduce to the case
where $\mu(M) \in \ZZ_{>0}$.
By Theorem~\ref{T:pure equiv}, there exists a pure $\phi$-module $M_0$
over $R^{\bd}$ with $M \cong M_0 \otimes R$.
By Proposition~\ref{P:h0}, we have $H^0(M_0) = H^0(M)$. 

Choose a basis $\be_1, \dots, \be_n$ of $M_0$ such that the matrix $A$
defined by
$\phi(\be_j) = \sum_i A_{ij} \be_i$ has all entries of valuation
at least $\mu(M)$. If $\bv = \sum c_i \be_i \in H^0(M)$ is nonzero, then
$c_i = \sum_j A_{ij} \phi(c_j)$ implies that $\min_i \{w(c_i)\} 
> \min_j\{w(c_j)\}$,
contradiction. Hence $H^0(M) = 0$.
\end{proof}
\begin{cor} \label{C:distinct slopes}
If $M$ and $N$ are pure $\phi$-modules over $R$ with
$\mu(M) < \mu(N)$, then $\Hom(M,N) = 0$.
\end{cor}
\begin{proof}
The conditions ensure that $M^\dual \otimes N$ is pure of positive slope;
by Lemma~\ref{L:h0 pos slope}, $\Hom(M,N) = H^0(M^\dual \otimes N) = 0$.
\end{proof}

\begin{theorem} \label{T:pure semi}
Let $M$ be a pure $\phi$-module over $R$ of slope $s$.
\begin{enumerate}
\item[(a)]
$M$ is module-semistable.
\item[(b)]
If $N$ is a $\phi$-submodule of $M$ with $\mu(N) = s$,
then $N$ is saturated, and both $N$ and $M/N$ are pure of slope $s$.
\end{enumerate}
\end{theorem}
\begin{proof}
For (a),
let $N$ be a $\phi$-submodule of $M$; we wish to show that $\mu(N) \geq s$.
By replacing $M$ by $\wedge^{\rank(N)} M$,
we may assume that $\rank(N) = 1$.
By Lemma~\ref{L:pure push}, we may assume further that $s \in \ZZ$.
By twisting, we may assume further that $N$ is trivial, so that $H^0(M) \neq 0$.
To avoid contradicting Lemma~\ref{L:h0 pos slope}, we must then have
$s \leq 0 = \mu(N)$, yielding semistability.

For (b), by applying $[a]_*$ and twisting,
we may again reduce to the case $s=0$. 
Let $M_0$ be an \'etale lattice in $M$; by Lemma~\ref{L:etale lattice},
the kernel of $M_0 \to M/N$ is a $\phi$-module over $R^{\inte}$,
so the image is as well. Let $P$ be the $R$-span of this
image; it is an \'etale $\phi$-submodule of $M/N$ of the same rank.
Since $\mu(N) = \mu(M) = 0$, we also have
$\mu(M/N) = 0$, so $M/N = P$ by Corollary~\ref{C:same rank}. Hence
$M/N$ is \'etale; the same logic applied after dualizing implies that
$N^\dual$ is \'etale, as then is $N$.
\end{proof}

\begin{cor} \label{C:split exact}
If $M_1,M_2$ are $\phi$-modules, then $M = M_1 \oplus M_2$ is pure of
slope $s$ if and only if both $M_1$ and $M_2$ are pure of slope $s$.
\end{cor}
\begin{proof}
If $M_1$ and $M_2$ are pure of the same slope, then visibly so is $M$.
Conversely, if $M$ is pure of slope $s$, then $M$ is semistable by
Theorem~\ref{T:pure semi}(a), so the $\phi$-submodules $M_1$ and $M_2$
each have slope at least $s$.
Since $\mu(M)$ is a weighted average of $\mu(M_1)$ and $\mu(M_2)$,
we must in fact have $\mu(M_1) = \mu(M_2) = s$; by 
Theorem~\ref{T:pure semi}(b), $M_1$ and $M_2$ are both pure of slope $s$.
\end{proof}

\begin{cor} \label{C:pure pull}
Let $M$ be a $\phi^a$-module over $R$. Then $M$ is pure of some slope $s$
if and only if $[a]^* M$ is pure of slope $s/a$.
\end{cor}
\begin{proof}
By Lemma~\ref{L:pure push}, $[a]^* M$ is pure of slope $s/a$ if and only if
$[a]_* [a]^* M$ is pure of slope $s$. If $M$ is pure of slope $s$,
then so are $(\phi^i)^* M$ for $i=0, \dots, a-1$; since
\begin{equation} \label{eq:push pull}
[a]_* [a]^* M \cong \oplus_{i=0}^{a-1} (\phi^i)^* M
\end{equation}
by Definition~\ref{D:pull push}, $[a]_* [a]^* M$ is pure of slope $s$.

Conversely, if $[a]_* [a]^* M$ is pure of slope $s$, then 
\eqref{eq:push pull} shows that $M$ is a direct summand of $[a]_* [a]^* M$,
and hence is pure by Corollary~\ref{C:split exact}.
\end{proof}

\subsection{The slope filtration theorem}

So far all of our work has been formal modulo the assumption of 
an appropriate analogue of Proposition~\ref{P:h1 both ways}.
We now restrict attention from general rings $R$ as in 
strong Hypothesis~\ref{hypo:robba}
to the Robba ring $\calR$ (as in Example~\ref{exa:Robba}), where one
can make the description of $\phi$-modules much more precise.

We have already described a natural  filtration on $\phi$-modules over
$\calR$, namely the Harder-Narasimhan filtration. The trouble is that
the construction is so formal that one cannot deduce any useful
properties about the resulting filtration or its associated slopes; for
instance, it is not clear that module-semistability is preserved by
tensor product. (The fact that the analogous statement is true for vector
bundles on smooth varieties in characteristic 0 is highly nontrivial:
it reduces to the case of tensoring two semistable
vector bundles of slope 0 on curves \cite{mr},
in which case it follows from an analytic classification of stable
bundles due to Narasimhan-Seshadri
\cite{ns1, ns2}.)
The slope filtration theorem, which is the main result of
this paper, asserts that in fact the
steps of the Harder-Narasimhan 
filtration are much more structured than one might have
otherwise predicted.
\begin{theorem}[Slope filtration theorem] \label{T:slope filt}
Every module-semistable $\phi$-module over the Robba ring $\calR$ is
pure. In particular, every $\phi$-module $M$ over $\calR$ admits a unique
filtration $0 = M_0 \subset M_1 \subset \cdots \subset M_l = M$ by
saturated $\phi$-modules whose successive quotients are pure with
$\mu(M_1/M_0)< \cdots < \mu(M_l/M_{l-1})$.
\end{theorem}
This theorem is stated as a forward reference, as its proof will occupy
most of the rest of the paper; here we give only a top-level summary.

\begin{proof}[Proof of Theorem~\ref{T:slope filt}]
The proof of Theorem~\ref{T:slope filt} will be obtained by constructing
(in Subsection~\ref{subsec:extended})
an extended Robba ring $\tcalR$ which also satisfies strong
Hypothesis~\ref{hypo:robba}, and then establishing the following facts.
\begin{itemize}
\item
If $M$ is a semistable $\phi$-module over $\calR$, then
$M \otimes \tcalR$ is also semistable (Theorem~\ref{T:ascend semistable}).
\item
If $\tilde{M}$ is a semistable $\phi$-module over $\tcalR$, then
$\tilde{M}$ is pure (Theorem~\ref{T:semistable pure}).
\item
If $M$ is a $\phi$-module over $\calR$ and
$M \otimes \tcalR$ is pure, then $M$ is pure
(Theorem~\ref{T:descend pure}).
\end{itemize}
These together yield the claim.
\end{proof}

\begin{remark}
Theorem~\ref{T:slope filt} implies that the tensor product of 
module-semistable $\phi$-modules is pure (by Corollary~\ref{C:tensor
pure}) and hence module-semistable (by Theorem~\ref{T:pure semi}). This
formally implies that the slopes of $\phi$-modules behave like 
valuations of eigenvalues, or like Deligne's weights in \'etale
cohomology. That is, if $M$ has slopes $c_1, \dots, c_m$ 
and $M'$ has slopes $c'_1, \dots, c'_n$, both counted with multiplicity,
then:
\begin{itemize}
\item
the slopes of $M \oplus M'$ are
$c_1, \dots, c_m, c'_1, \dots, c'_n$;
\item
the slopes of $M \otimes M'$ are
$c_i c'_j$ for $i=1, \dots, m$ and $j = 1, \dots, n$;
\item
the slopes of $\wedge^d M$ are
$c_{i_1} + \cdots + c_{i_d}$ for $1 \leq i_1 < \cdots < i_d \leq m$;
\item
the slopes of $[a]_* M$ are $ac_1, \dots, ac_m$;
\item
the slopes of $M(b)$ are $c_1 + b, \dots, c_m + b$;
\item
the slopes of $[a]^* M$ are $c_1/a, \dots, c_m/a$, each repeated $a$
times.
\end{itemize}
In some sense, the slope filtration theorem is thus playing a role in this
theory analogous to Deligne's analysis of determinantal weights in his
second proof of the Weil conjectures \cite{deligne}.
\end{remark}

\begin{remark}
The uniqueness in Theorem~\ref{T:slope filt} means that the slope filtration
inherits any additional group action on the original $\phi$-module. In 
particular, if $M$ is a $(\phi, \Gamma)$-module, then the steps of the
slope filtration are $(\phi, \Gamma)$-submodules of $M$. As shown by Berger
\cite[Th\'eor\`eme~V.2.1]{berger-weak},
this leads to a
proof of the Colmez-Fontaine theorem that $(\phi, N)$-modules over a $p$-adic
field which are \emph{weakly admissible}, in the sense of satisfying a 
necessary numerical criterion, indeed arise
from Galois representations via $p$-adic Hodge theory.
(See also the variant of Berger's argument given by Kisin \cite{kisin}.)
\end{remark}

\begin{remark} \label{R:cherbonnier}
The \'etale $(\phi, \Gamma)$-modules attached to Galois representations of a
$p$-adic field were originally defined by Fontaine over the $p$-adic completion
of $\calR^{\bd}$; the fact that they can be descended to
$\calR^{\bd}$ is a theorem of Cherbonnier and Colmez
\cite[Corollaire~III.5.2]{cherbonnier-colmez}.
The fact that the descent is unique follows from the fact that
the base change from \'etale $\phi$-modules over $\calR^{\bd}$ to its 
completion is fully faithful, which in turn follows from
Proposition~\ref{P:cherbonnier1}.
\end{remark}

\begin{remark}
In the context of $p$-adic differential equations and rigid cohomology,
Theorem~\ref{T:slope filt} arises with $M$ carrying the extra structure of
a connection $\nabla: M \to M \otimes \Omega^1_{\calR/K}$ compatible with 
the $\phi$-action; that is, $M$ is a \emph{($\phi, \nabla$)-module}.
One can see that the steps of the slope filtration are
$(\phi, \nabla)$-submodules by using Corollary~\ref{C:distinct
slopes} as follows.
The map $\nabla$ induces a homomorphism $M_1 \to (M/M_1) \otimes 
\Omega^1_{\calR/K}$
of $\phi$-modules. Since $\Omega^1_{\calR/K}$ 
is a rank $1$ $\phi$-module of nonnegative
slope (the slope is actually positive, but we don't need this here), each slope
of $(M/M_1) \otimes \Omega^1_{\calR/K}$ is strictly greater than $\mu(M_1)$. Repeated
application of Corollary~\ref{C:distinct slopes} yields the claim.

Given that the slope filtration is a filtration by $(\phi,\nabla)$-submodules,
one may prove the local monodromy theorem for $p$-adic differential equations
as in \cite{kedlaya-local},
by showing each successive quotient in the slope filtration becomes trivial
as a $\nabla$-module after tensoring with a suitable finite unramified extension
of $\calR^{\inte}$. This reduces easily to the \'etale case, which is
a theorem of Tsuzuki \cite[Theorem~4.2.6]{tsuzuki}. Beware, however, that this
last step only applies for $\phi_K$ absolute; in particular, this approach
cannot be used to prove \cite[Proposition~6.2.1]{berger-colmez}.
\end{remark}

\begin{remark} \label{R:hartl-pink}
By \cite[Theorem~11.1]{hartl-pink},
the conclusion of Theorem~\ref{T:slope filt} also holds
in the situation of Example~\ref{ex:hartl-pink};
indeed, what one obtains is an analogue of the classification
of $\phi$-modules over the extended Robba ring $\tcalR$ to be introduced
in Section~\ref{sec:dm}.
That result is not covered by this paper, though (as \cite{hartl-pink} 
already points out) there are very strong parallels between the ensuing
calculations. However, Theorem~\ref{T:slope filt} itself
does address a related
situation: if we take $K = k((z))$ with $k$ of characteristic $p>0$, 
and $\phi_K$ to be a power of the absolute Frobenius,
then $\calR$ consists of Laurent series
in $t$ over $z$ which converge for $|z|^c < |t| < 1$ for some $c>0$.
Since the valuation on $k$ is trivial, it is equivalent to require
convergence when $0 < |z| < |t|^{1/c}$; that is, we are considering series in 
$z$ over $k((t))$ convergent on some punctured disc around the origin.
In this case (assuming $q$ is a power of $p$), Theorem~\ref{T:slope filt}
is a result of Hartl \cite[Theorem~1.7.7]{hartl}.
\end{remark}

It would be interesting to know about the following $q$-analogue of 
Remark~\ref{R:hartl-pink}; it may be related to the formal classification
of linear difference operators \cite{praagman}, in much the same way that
the construction of the canonical lattice of an irregular meromorphic
connection
\cite{malgrange} reduces to the formal classification of linear
differential operators \cite{levelt}.

\begin{question} \label{Q:q-analogue}
Let $K$ be a complete field, either archimedean or nonarchimedean.
Take $R$ to be the ring of germs of 
analytic functions over $K$ on punctured discs
around the origin, $R^{\bd}$ to be the germs meromorphic at zero,
$w$ to be the order of vanishing at zero, and
$\phi$ to be the map $\sum c_i t^i \mapsto \sum c_i q^i t^i$ for some $q \in 
K^*$ with $|q| < 1$.
Does the analogue of Theorem~\ref{T:slope filt} hold in this setting?
\end{question}

\begin{remark} \label{R:dieudonne-manin}
The conclusion of Theorem~\ref{T:slope filt} also holds for
$\phi$-modules over $K$ itself; this is a straightforward consequence of
Proposition~\ref{P:factor twisted}.
In addition, if $\phi$ is bijective on $K$, then it is easy to check that
$H^1(M) = 0$ for $M$ pure of nonzero slope, so the slope filtration
splits uniquely. This gives a semilinear
 analogue of the eigenspace decomposition
of a vector space equipped with a linear transformation.
If $k$ is algebraically closed of 
characteristic $p>0$ and $\phi$ is an absolute Frobenius lift,
this recovers the Dieudonn\'e-Manin classification of rational Dieudonn\'e
modules \cite{manin}.
\end{remark}

\begin{remark} \label{R:no identity}
The conclusion of Theorem~\ref{T:slope filt} does not hold for $\phi$
equal to the identity map on $\calR$. In fact, the conclusion is equivalent
to the condition that the characteristic polynomial of $\phi$ have all
coefficients in $\calR^{\bd}$, whereas the definition of a $\phi$-module
only forces the determinant to belong to $\calR^{\bd}$.
\end{remark}

\section{Classification over an extended Robba ring}
\label{sec:dm}

In this section and the next, we give a proof of
Theorem~\ref{T:slope filt}. Although somewhat simplified in some technical
aspects, 
the argument follows the same arc as in
\cite{kedlaya-local} and \cite{kedlaya-slope}, with two basic stages.
In the first stage, performed in this section, we show that $\phi$-modules
over a suitable overring of $\calR$ admit a very simple classification
(analogous to the Dieudonn\'e-Manin classification alluded to
in Remark~\ref{R:dieudonne-manin}),
and in particular admit a slope filtration.
In the second stage, we show that the slope filtration descends back to
$\calR$.

On a first reading, we recommend reading only Subsection~\ref{subsec:overview}
for an overview, then returning later for the technical details
in the rest of the section.

\subsection{Overview}
\label{subsec:overview}

\begin{hypo} \label{H:field size}
Throughout this section, assume that $\phi$ is a relative Frobenius lift
on $\calR$ such that $\phi_K$ is an \emph{automorphism}
of $K$. Also assume that 
any \'etale $\phi$-module over $K$ is trivial; this is equivalent
to asking that any $\phi$-module over the residue field
$k$ be trivial. It also implies
that $H^1$ vanishes for any \'etale $\phi$-module over $K$ or any
$\phi$-module over $k$.
Using Definition~\ref{D:h0 h1}, 
we  deduce the same conclusions with $\phi$ replaced by
$\phi^a$ for any positive integer $a$.
\end{hypo}

\begin{remark}
In the absolute Frobenius case, Hypothesis~\ref{H:field size}
can be satisfied by taking $k$ to be algebraically closed.
In general,
one must work a bit harder; see Proposition~\ref{P:splitting field}.
\end{remark}

We will define (Definition~\ref{D:extended Robba})
an \emph{extended Robba ring} $\tcalR$ which has
the following properties:
\begin{itemize}
\item
$\tcalR$ is a B\'ezout domain containing $\calR$, and admits an automorphism
$\phi$ extending the given Frobenius lift on $\calR$
(see Remark~\ref{R:analytic ring} and Proposition~\ref{P:embedding}).
\item
The units in $\tcalR$ are the nonzero elements of a subfield
$\tcalR^{\bd}$, which is the fraction field of a discrete valuation
ring $\tcalR^{\inte}$ for which $\tcalR^{\inte} \cap \calR = \calR^{\inte}$
(see Remark~\ref{R:analytic ring}).
\item
The strong form of Hypothesis~\ref{hypo:robba} holds for $R = \tcalR$
(see Proposition~\ref{P:h1 both ways2}).
\end{itemize}

The classification of $\phi$-modules over $\tcalR$
rests on a sequence of structural results, which we state in roughly
increasing order of difficulty; their proofs occupy the remainder of this
section.

\begin{prop} \label{P:h0 nonzero}
Let $M,N$ be pure $\phi$-modules over $\tcalR$ obtained by base change
from $K$, with $\mu(M) > \mu(N)$. Then $\Hom(M,N) \neq 0$.
\end{prop}
\begin{proof}
See Subsection~\ref{subsec:extended}.
\end{proof}

\begin{notation}
Choose a uniformizer $\pi$ of $K$, and let $\tcalR(1)$ be the
$\phi$-module of rank 1 and degree 1 on which $\phi$ acts on some
generator via multiplication by $\pi$. We use $\tcalR(1)$ as a
twisting sheaf, writing
$M(n) = M \otimes \tcalR(1)^{\otimes n}$.
\end{notation}

\begin{prop} \label{P:eigenvector}
Let $M$ be a nonzero $\phi$-module over $\tcalR$. Then for all sufficiently
large integers $n$, $H^0(M(-n)) \neq 0$ and $H^1(M(-n)) = 0$.
\end{prop}
\begin{proof}
See Subsection~\ref{subsec:fixed vec}.
\end{proof}

\begin{prop} \label{P:pure equiv2}
For any rational number $s$, the base change functor from pure
$\phi$-modules of slope $s$ over $K$ to pure $\phi$-modules of
slope $s$ over  $\tcalR$ is an equivalence of categories.
\end{prop}
\begin{proof}
See Subsection~\ref{subsec:classif}.
\end{proof}

\begin{prop} \label{P:local calc}
Let $n$ be a positive integer,
let $N'$ be a pure $\phi^n$-module over $\tcalR$ of rank $1$ and degree $1$,
let $P$ be a pure $\phi$-module over $\tcalR$ of rank $1$ and degree $-1$,
and suppose
\[
0 \to [n]^* N' \to M \to P \to 0
\]
is a short exact sequence of $\phi$-modules. Then
$H^0(M) \neq 0$.
\end{prop}
\begin{proof}
See Subsection~\ref{subsec:local calc}.
\end{proof}

These assemble to give the following classification theorem.
\begin{theorem} \label{T:semistable pure}
Any module-semistable $\phi$-module over $\tcalR$ is pure.
Consequently, the
successive quotients of the HN filtration of a $\phi$-module over $\tcalR$
are all pure.
\end{theorem}
\begin{remark} \label{R:any hom}
Before proving Theorem~\ref{T:semistable pure}, we make an observation
which figures prominently in the argument.
If one knows Theorem~\ref{T:semistable pure} for $\phi$-modules
of rank $\leq n$, it follows from 
Propositions~\ref{P:h0 nonzero} and~\ref{P:pure equiv2} 
(and the assumption that \'etale $\phi$-modules over $K$ are trivial,
as built into Hypothesis~\ref{H:field size})
that
for $M$ a pure $\phi$-module over $\tcalR$ and $N$ an \emph{arbitrary}
$\phi$-module over $\tcalR$ with $\rank(N) \leq n$ and 
$\mu(M) \geq \mu(N)$, we have
$\Hom(M,N) \neq 0$; in particular, if $\rank(M) = 1$, we would have
an injection of $M$ into $N$. This is because the first step of the HN
filtration of $N$ always has slope $\leq \mu(N)$.
\end{remark}

\begin{proof}[Proof of Theorem~\ref{T:semistable pure}]
We proceed by induction on rank, the case of rank 1 being evident.
Assume that $n \geq 1$ and that for every positive integer $a$,
every semistable $\phi^a$-module of rank $\leq n$
is pure.  
Suppose that $M$ is a semistable $\phi^a$-module of rank $n+1$ over $\tcalR$;
we wish to show that $M$ is pure. 
We may reduce to the case where
$\mu(M) \in \ZZ$ by applying $[d]_*$ and invoking
Corollary~\ref{C:semi push} (to see that semistability is preserved)
and Lemma~\ref{L:pure push} (to see that purity is reflected);
we may then twist to ensure $\mu(M) = 0$.
For ease of notation, we will assume hereafter that $M$ is a $\phi$-module
(at the expense of replacing $\phi$ by a power, which does not
disturb Hypothesis~\ref{H:field size}).

Put $M' = [n]_* M$; then $M'$ is semistable by Corollary~\ref{C:semi
push} again. 
By Proposition~\ref{P:eigenvector}, there exists a nonnegative integer $c$
such that $M'$ admits a pure $\phi^{n}$-submodule $N'$ 
of rank 1 and slope $c$; choose $c$ as small as possible.
Suppose that $c \geq 2$; since $\mu(M'/N') < 0 \leq c-2$, we may apply
Remark~\ref{R:any hom} to produce a 
$\phi^{n}$-submodule of $M'/N'$ isomorphic to 
$\tcalR(c-2)$. Let $Q'$ be the inverse image of that submodule in $M'$;
applying Proposition~\ref{P:local calc} (in the case $n=1$) 
to the exact sequence
\[
0 \to N'(1-c) \to Q'(1-c) \to \tcalR(-1) \to 0,
\]
we see that $H^0(Q'(1-c)) \neq 0$ and hence $H^0(M'(1-c)) \neq 0$, contradicting
the minimality of $c$.

Suppose that $c = 1$. Put $N = [n]^* N'$; 
then $N$ is pure of slope $1/ n$ by Corollary~\ref{C:pure pull}.
The adjunction between $[n]^*$ and $[n]_*$
converts the inclusion $N' \hookrightarrow M'$ into a nonzero map
$f: N \to M$. Since $N$ is semistable by Theorem~\ref{T:pure semi},
$\mu(f(N)) \leq 1/n$;
moreover, the denominator of $\mu(f(N))$ is at most $\rank(f(N)) \leq n$. 
Consequently, either $\mu(f(N)) \leq 0$, in which case
Remark~\ref{R:any hom} implies that $H^0(f(N)) \neq 0$; or $\mu(f(N)) =
1/n$,
in which case $f$ must be injective and we have an exact sequence
\[
0 \to N \to M \to P \to 0
\]
with $P$ pure of rank 1 and slope $-1$, to which we apply
Proposition~\ref{P:local calc} to deduce that $H^0(M) \neq 0$.
In either case, we contradict the minimality of $c$.

We deduce that $c=0$, i.e., $M'$ admits a nontrivial \'etale $\phi$-submodule 
$N'$; the quotient $M'/N'$ is also semistable, hence pure by
the induction hypothesis. By Theorem~\ref{T:pure ext}, $M'$ is pure, as then
is $M$ by Lemma~\ref{L:pure push}. This completes the proof.
\end{proof}

\begin{remark}
In the proof of Theorem~\ref{T:semistable pure},
the passage from $M$ to $M'$ is made in order to simplify the statement
of Proposition~\ref{P:local calc}. One can do some extra work 
to prove a version of Proposition~\ref{P:local calc} in which $[n]^* N$
is replaced by any pure $\phi$-module of rank $n$ 
and degree 1; however, the internal
improvement is immaterial in the end, as even this stronger
form of Proposition~\ref{P:local calc}
is itself an immediate consequence of Theorem~\ref{T:semistable pure}.
\end{remark}

\subsection{The extended Robba ring}
\label{subsec:extended}

We now go back and construct the extended Robba ring $\tcalR$.

\begin{defn}
Let $R$ be a ring and let $G$ be a totally ordered abelian group. The
ring of \emph{Hahn series} (or \emph{Mal'cev-Neumann series}, or
\emph{generalized power series}) over $R$ with value group $G$ is the
set of functions $f: G \to R$ with well-ordered support, with pointwise
addition and multiplication
given by convolution; it is a standard calculation 
\cite[Chapter~13]{passman} to verify that these operations give
a well-defined ring, which is a field if $R$ is.
We typically represent elements of this ring as formal
series $\sum_{g \in G} r_g u^g$ in some dummy variable $u$ with powers
indexed by $g \in G$, and the ring is correspondingly denoted $R((u^G))$.
For $G \subseteq \RR$,
we view $R((u^G))$ as being equipped with the $u$-adic valuation $v$ sending
$\sum_g r_g u^g$ to the smallest $g$ for which $r_g \neq 0$ (i.e., the
least element of the support).
\end{defn}

\begin{lemma} \label{L:kernel cokernel}
Let $\phi: R((u^\QQ)) \to R((u^\QQ))$ be an automorphism of the form
$\sum_i a_i u^i \mapsto \sum_i \phi_R(a_i) 
u^{qi}$, with $\phi_R$ an automorphism
 of $R$. Then the map $1 - \phi$ is bijective on the set of series with
 zero constant term.
\end{lemma}
\begin{proof}
If $x \in R((u^\QQ))$ and $v(x) < 0$, then
$v(x - \phi(x)) = q v(x)$, whereas if $v(x)> 0$, then $v(x - 
\phi(x)) = v(x)$.
This proves injectivity.

Given $x \in R((u^\QQ))$, write $x = \sum_i x_i u^i$, and put
\begin{align*}
y_+ &= \sum_{j=0}^\infty \sum_{i>0} \phi^j_R(x_i) u^{iq^j} \\
y_- &= \sum_{i<0} \left( \sum_{j=0}^\infty 
-\phi^{-j-1}_R(x_{iq^{j+1}}) \right) u^{i}.
\end{align*}
Since both sums give well-defined elements of $R((u^\QQ))$ (in the definition of $y_-$, the
sum over $j$ is finite for each $i$), we may
put $y = y_+ + y_-$, which has zero constant term and satisfies
$y - \phi(y) = x - x_0$. This proves surjectivity.
\end{proof}

\begin{cor} \label{C:kernel cokernel}
With $k$ as in Hypothesis~\ref{H:field size},
for any $c \in k^*$, the map $1 - c \phi$ on $k((u^\QQ))$ is surjective.
\end{cor}
\begin{proof}
By Hypothesis~\ref{H:field size}, there exists $a \in k^*$ such that
$\phi(a) = c a$, so we can always write 
\[
(1 - c\phi)(x) = a^{-1} (ax - \phi(ax)).
\]
It thus suffices to check the case $c=1$; this follows from
Lemma~\ref{L:kernel cokernel} and the fact that $1-\phi$ is surjective
on $k$, which again is a consequence of Hypothesis~\ref{H:field size}.
\end{proof}

Corresponding to the extension from power series to generalized power series,
we define an enlargement of the Robba ring. We first construct the ring,
then the embedding of the original Robba ring into it.

\begin{defn} \label{D:extended Robba}
For $r>0$, 
let $\tcalR^r$ be the set of formal sums $\sum_{i \in \QQ} a_i u^i$
with $a_i \in K$, satisfying the following conditions.
\begin{itemize}
\item
For each $c > 0$, the set of $i \in \QQ$ such that $|a_i| \geq c$ is well-ordered.
\item
We have $|a_i| e^{-ri} \to 0$ as $i \to -\infty$.
\item
For all $s>0$, we have $|a_i| e^{-si} \to 0$ as $i \to +\infty$.
\end{itemize}
Then $\tcalR^r$ can be shown to form a ring.
We call the union $\tcalR = \tcalR_K =  \cup_{r} \tcalR^r$ the
\emph{extended Robba ring} over $K$. 
Let  $\tcalR^{\bd}$ and $\tcalR^{\inte}$ be the subrings
 of $\tcalR$ consisting of series
with bounded and integral coefficients, respectively. 
We equip $\tcalR^r$ with the norm
\[
\left| \sum_i a_i u^i \right|_r = \sup_i \{|a_i| e^{-ri}\}
\]
and $\tcalR$ with the automorphism
\[
\phi \left(\sum_i a_i u^i \right) = \sum_i \phi_K(a_i) u^{qi}.
\]
\end{defn}

\begin{remark} \label{R:analytic ring}
The ring $\tcalR$
can be viewed as an example of an ``analytic ring'' in the sense of
\cite[\S 2.4]{kedlaya-slope}, by taking $\phi_K$ to be an absolute Frobenius lift
on $K$.
Thus the results of \cite[Chapter~2]{kedlaya-slope} apply to show that
$\tcalR$ shares many of the nice properties of $\calR$, as follows.
\begin{itemize}
\item
The ring $\tcalR$ is a B\'ezout domain \cite[Theorem~2.9.6]{kedlaya-slope}.
\item
The ring $\tcalR^{\inte}$ is 
a henselian discrete valuation ring, and its fraction field
is $\tcalR^{\bd}$
\cite[Lemma~2.1.12]{kedlaya-slope}.
\item
The units of $\tcalR$ are the nonzero elements of $\tcalR^{\bd}$
\cite[Lemma~2.4.7]{kedlaya-slope}.
\end{itemize}
\end{remark}

\begin{prop} \label{P:embedding}
There exists a $\phi$-equivariant embedding $\psi: \calR \hookrightarrow
\tcalR$ such that for any $r_0$ as in Remark~\ref{R:annuli} and any
$r \in (0,r_0)$, $\calR^r$ maps to $\tcalR^r$ preserving $|\cdot|_r$.
\end{prop}
\begin{proof}
We inductively construct homomorphisms $\psi_l: \calR \to \tcalR$,
each of the form $\psi_l(\sum c_i t^i) = \sum c_i u_l^i$ for some
$u_l \in \tcalR^{\inte}$ with $|u_l|_r = |t|_r$ for $r \in (0,r_0)$, satisfying
\[
\psi_l(\phi(x)) \equiv \phi(\psi_l(x)) \pmod{\pi^l} 
\qquad (x \in \calR^{\inte}),
\]
starting with $u_1 = u$.
Given $\psi_l$, we may repeatedly invoke
Corollary~\ref{C:kernel cokernel} (if $q \neq 0$ in $k$)
or the fact that $\phi$ is surjective on $\tcalR^{\inte}$ 
(if $q =0$ in $k$)
to construct $\Delta \in \tcalR^{\inte}$ with
\begin{equation} \label{eq:Delta}
\phi(\pi^l \Delta/u ) - q (\pi^l \Delta/u) 
= (\psi_l(\phi(t)) - \phi(u_l))/u^q.
\end{equation}
For any $r \in (0,r_0)$,
\[
|\psi_l(\phi(t))|_{r/q}, |\phi(u_l)|_{r/q} \leq |t^q|_{r} = |u^q|_{r}
\]
and so the right side of \eqref{eq:Delta} has $(r/q)$-norm at most 1.
{}From this plus either
the proof of Lemma~\ref{L:kernel cokernel} (if $q \neq 0$ in $k$)
or direct inspection (if $q=0$ in $k$),
we deduce that $|\pi^l \Delta/u|_{r} \leq 1$. 
We may thus set $u_{l+1} = u_l + \pi^l \Delta$ to construct $\psi_{l+1}$;
this has the desired effect because
\[
\psi_{l+1}(\phi(t))
\equiv \psi_l(\phi(t)) + q \pi^l \Delta u^{q-1} \pmod{\pi^{l+1}}.
\]

The property $|u_l|_r = |t|_r$ implies that
each $\psi_l$ carries
$\calR^r$ to $\tcalR^r$ preserving $|\cdot|_r$.
By continuity, we obtain a map $\psi$ with the same property,
as desired.
\end{proof}

\begin{lemma} \label{L:where fixed}
The fixed elements of $\tcalR$ under $\phi$ all belong to $K$.
\end{lemma}
\begin{proof}
For $x =\sum_i a_i u^i \in\tcalR$, we have $\phi(x) = \sum_i \phi_K(a_i) u^{qi}$.
If $\phi(x) = x$ and $a_i \neq 0$ for some $i \neq 0$, then $|a_{i q^n}| = |a_i|$
for all $n \in \ZZ$; but this contradicts the fact that for any $c>0$, the
set of $i \in\QQ$ with $|a_i| \geq c$ is well-ordered.
Hence $a_i = 0$ for all $i \neq 0$, proving the claim.
\end{proof}

We now notice that strong Hypothesis~\ref{hypo:robba} holds for $\tcalR$.
\begin{prop} \label{P:h1 both ways2}
Let $A$ be an $n \times n$ matrix over $\tcalR^{\inte}$. Then the map
$\bv \mapsto \bv - A \phi(\bv)$ on column vectors induces a bijection
on $(\tcalR/\tcalR^{\bd})^n$.
\end{prop}
\begin{proof}
The proof proceeds as in Proposition~\ref{P:h1 both ways}, using the
definition of $|\cdot|_r$ given in Definition~\ref{D:extended Robba}.
\end{proof}

\begin{remark}
As a reminder, here are some key properties of
$\tcalR$ which we will use going forward.
\begin{itemize}
\item
Given a relative Frobenius lift $\phi$ on $\calR$, we can define an
action of $\phi$ on $\tcalR$ and an 
equivariant embedding $\psi: \calR \hookrightarrow
\tcalR$ which preserves $|\cdot|_r$ for $r \in (0,r_0)$ 
(Proposition~\ref{P:embedding}).
\item
The map $\phi$ is bijective on $\tcalR$.
\item
The map $1 - \phi$ is bijective on
$\tcalR^{\inte}/ \gotho_K$ (easy consequence of
Lemma~\ref{L:kernel cokernel}).
\item
There is a natural direct limit topology, restricting to the direct limit
of Fr\'echet topologies on $\calR$, under which $\tcalR$ is complete.
\end{itemize}
In \cite{kedlaya-local} and \cite{kedlaya-slope}, the role of our $\tilde{\calR}$
is played by the ring $\Gamma^{\alg}_{\an,\con}$, which is constructed
to be minimal for the
above properties; that ring 
coincides with the ring 
denoted $\tilde{\mathbf{B}}^\dagger_{\mathrm{an}}$
(as in \cite[\S II]{berger-cst}) or
more commonly $\tilde{\mathbf{B}}^\dagger_{\mathrm{rig}}$
(as in \cite{colmez-bourbaki}).
We opt here for the ring $\tcalR$ instead in hopes that
the construction using generalized power series makes the analogy to
$\calR$ a bit more apparent.
\end{remark}

To conclude this section, we prove  Proposition~\ref{P:h0 nonzero}:
if $M,N$ are pure $\phi$-modules over $\tcalR$ obtained by base
change from $K$, with $\mu(M) > \mu(N)$, then $\Hom(M,N) \neq 0$.
\begin{proof}[Proof of Proposition~\ref{P:h0 nonzero}]
It is equivalent to show that if $M$ is pure with $\mu(M) < 0$, obtained
by base change from $K$, then
$H^0(M) \neq 0$. Write 
$M = M_0 \otimes_K \tcalR$ for $M_0$ a pure $\phi$-module over $K$.
Take any nonzero $\bw \in M_0$ and any $i>0$;
the sum
\[
\bv = \sum_{n \in \ZZ} \phi^n (u^i \bw)
\]
will converge to a nonzero element of $H^0(M)$.
(Compare \cite[Proposition~3.3.4(c2)]{kedlaya-slope}.)
\end{proof}

\subsection{Construction of fixed vectors}
\label{subsec:fixed vec}

We next treat Proposition~\ref{P:eigenvector}:
if $M$ is a nonzero $\phi$-module over $\tcalR$, then for all sufficiently
large integers $n$, $H^0(M(-n)) \neq 0$ and $H^1(M(-n)) = 0$.
(Also compare \cite[Theorem~4.1]{hartl-pink}.)
\begin{proof}[Proof of Proposition~\ref{P:eigenvector}]
We follow \cite[Proposition~4.2.2]{kedlaya-slope}.
View $M$ as a space of column vectors with the action of $\phi$ given by
multiplication by the matrix $A$ times the componentwise action.
Choose $r>0$ so that $A$ and $A^{-1}$ have entries in $\tcalR^{qr}$.

For $d \in \QQ_{>0}$ to be specified below, 
define the ``splitting functions''
$f_d^+, f_d^-$ as follows: given $x = \sum a_i u^i$, put
\[
f_d^+(x) = \sum_{i \geq d} a_i u^i, \qquad
f_d^-(x) = \sum_{i<d} a_i u^i,
\]
then extend to vectors componentwise. For $\bw$ a vector, we write
$\bw^{\pm}$ for $f_d^{\pm}(\bw)$.

Define the map $g: M \to M$ by
\[
g(\bw) = \pi^{-n} A \phi(\bw^+) + \phi^{-1}(\pi^n A^{-1} \bw^-)
\]
and note that
\begin{equation} \label{eq:bound g}
|g(\bw)|_r \leq \max\{|\pi|^{-n} |A|_r e^{-rd(q-1)},
|\pi|^{n} |A^{-1}|_{qr} e^{-rd(q^{-1}-1)} \} |\bw|_r.
\end{equation}
If we can choose $d$ such that
the two quantities in the maximum in \eqref{eq:bound g} are both
strictly less than 1, then $g$ will be contractive towards zero. This happens
if
\begin{equation} \label{eq:bound d}
d \in \left( \frac{n(-\log |\pi|) + \log |A|_r}{r(q-1)},
\frac{qn(-\log |\pi|) - q \log |A^{-1}|_{qr}}{r(q-1)} \right);
\end{equation}
for $n$ sufficiently large the interval is nonempty.
(Note  that
consistently with Proposition~\ref{P:h0 nonzero}, if $M$ is \'etale over $K$
we can take any $n>0$.)

Fix $n,d$ satisfying \eqref{eq:bound d}.
Given $\bw$ with entries in $\tcalR^r$, we define the sequence
$\bw_0 = \bw, \bw_{l+1} = g(\bw_l)$, then set
\begin{equation} \label{eq:sum}
\bv = \sum_{l=0}^\infty (\bw_l^+ - \phi^{-1} (\pi^n A^{-1} \bw_l^-)),
\end{equation}
so that $|\bv|_r \leq |\bw|_r$ and $\bv - \pi^{-n} A \phi(\bv) = \bw$.
We only know \emph{a priori} that the sum defining $\bv$ converges under
$|\cdot|_r$, but using 
the equation $\bv = \pi^{-n} A \phi(\bv) + \bw$, we may deduce that
the sum converges under $|\cdot|_{r/q}$, $|\cdot|_{r/q^2}$, and so on.
Hence $\bv$ has entries in $\tcalR^r$, yielding
$H^1(M(-n)) = 0$.

To deduce $H^0(M(-n)) \neq 0$,
we modify the previous construction slightly.
Put $\bw = (u^d, 0, \dots, 0)$ and construct $\bv$ as in
\eqref{eq:sum}.
Then put $\bw'_0 = \bw$, $\bw'_1 = 
\phi^{-1}(\pi^n A^{-1} \bw'_0)$, and 
$\bw'_{l+1} = g(\bw'_l)$ for $l \geq 1$. 
(That is, at the first step, transfer
the boundary term $u^d$ from the plus part to the minus part.)
If we now define
\[
\bv' = -\phi^{-1}(\pi^n A^{-1}\bw'_0) + 
\sum_{l=1}^\infty ((\bw_l')^+ - \phi^{-1} (\pi^n A^{-1} (\bw_l')^-)),
\]
we obtain $\bv' - \pi^{-n} A \phi(\bv') = \bw$ as before. However,
$|\bv|_r = |u^d|_r$ whereas $|\bv'|_r < |u^d|_r$, so 
$\bv- \bv'$ is a nonzero element of $H^0(M(-n))$, as desired.
\end{proof}

\subsection{Twisted polynomials and their Newton polygons}

Before continuing, we need to analogize, to the realm of twisted polynomials
over $k((u^\QQ))$, some facts about polynomials over
valued fields and their Newton polygons. With a bit of care, we can obtain
at the same time some results over $K$ which we will need later
(see Proposition~\ref{P:splitting field}).

\begin{notation} \label{N:valued field}
Throughout this subsection only, fix a real number $s \geq 1$, 
and
let $F$ be a field equipped with an automorphism $\phi = \phi_F$
and a valuation $v_F$ with the properties that $F$ is complete under $v_F$
and $v_F(\phi_F(x)) = s v_F(x)$ for all $x \in F$.
Let $\gotho_F$ and $\gothm_F$
denote the valuation subring of $F$ and the maximal ideal of $\gotho_F$, 
respectively.
\end{notation}

\begin{defn} \label{D:Newton poly}
For $i \in \ZZ$, write $[i] = \sum_{j=0}^{i-1} s^j$, so that 
$[0] = 0$, $[1] = 1$, and $[i+j] = [i] + s^i[j]$. 
For $r \in \RR$ and $P(T) \in F\{T^{\pm}\}$,
write $P(T) = \sum_{i \in \ZZ} a_i T^i$, and write
\[
v_r(P) = \min_i \{ v_F(a_i) + r[i]\}.
\]
Define the \emph{homogeneous Newton polygon} of $P$ as the lower convex
hull of the set
\[
\{(-[i], v_F(a_i)): i \in \ZZ\};
\]
we refer to the slopes of this polygon as the \emph{(Newton) slopes of $P$}.
\end{defn}

\begin{lemma} \label{L:scale product}
For $P(T) \in F\{T\}$ and $Q(T) \in F\{T^{-1}\}$
such that $v_r(Q) \geq 0$, we have
$v_r(PQ) \geq v_r(P) + v_r(Q)$.
\end{lemma}
\begin{proof}
Write $P(T) = \sum_{i \geq 0} a_i T^i$ and $Q(T) = 
\sum_{j \leq 0} b_j T^j$.
We have
\[
(PQ)(T) = \sum_k \sum_{i+j=k} a_i \phi^i(b_j) T^k,
\]
and
\begin{equation} \label{eq:newton}
v_F(a_i \phi^i(b_j)) + [i+j]r =
v_F(a_i) + [i]r + s^i (v_F(b_j) + [j]r).
\end{equation}
The right side of \eqref{eq:newton} is at least
$v_r(P) + s^i v_r(Q)$. Since $i \geq 0$ and $s \geq 1$, if $v_r(Q) \geq 0$, then 
the right side of \eqref{eq:newton} is at least $v_{r}(P) + v_{r}(Q)$.
This yields the claim.
\end{proof}

\begin{prop} 
Let $r_0 \in \RR$ be a real number, and suppose that
$P(T) \in F\{T\}$ and $Q(T) \in F\{T^{-1}\}$ are such that
$P$ has constant term $1$ and all slopes $\leq r_0$,
and $Q$ has constant term $1$ and all slopes $\geq r_0$.
Then
the slopes of $PQ$ are obtained by taking the union (with multiplicities)
of the sets of slopes
of $P$ and $Q$.
\end{prop}
\begin{proof}
The conditions on the slopes of $P$ and $Q$ imply that
\begin{align*}
r \geq r_0 &\implies v_r(P) = 0, v_r(Q) \leq 0 \\
r \leq r_0 &\implies v_r(P) \leq 0, v_r(Q) = 0.
\end{align*}
It thus suffices to check that
\[
v_{r}(PQ)  = \begin{cases} v_{r}(Q) & r > r_0 \\
0 & r = r_0 \\
v_{r}(P) & r < r_0.
	 \end{cases}
\]
Retain notation as in Lemma~\ref{L:scale product}.
If $r \geq r_0$, take the smallest $j$ that minimizes
$v_F(b_j) + [j]r$; then \eqref{eq:newton} equals $v_{r}(Q)$ for $i=0$ but not
 for any other pair $i,j$ with the same sum. 
If $r \leq r_0$, take the largest $i$ that
minimizes $v_F(a_i) + [i]r$; then \eqref{eq:newton} equals $v_{r}(P)$ for $j=0$
but not for any other pair $i,j$ with the same sum.
This yields the desired result.
\end{proof}

\begin{prop} \label{P:factor twisted}
Let $r \in \RR$ be a real number, and suppose that
$R \in F\{T^{\pm}\}$
satisfies $v_r(R-1) > 0$.
Then there exist $c \in F$, $P(T) \in F\{T\}$,
$Q(T) \in F\{T^{-1}\}$ such that
$v_F(c-1) > 0$,
$P$ has constant term $1$ and all slopes $< r$,
$Q$ has constant term $1$ and all slopes $> r$,
and $cPQ = R$.
\end{prop}
\begin{proof}
Put $c_0 = P_0 = Q_0 = 1$. Given $c_i, P_i, Q_i$, write
$R - c_i P_i Q_i = \sum_j r_j T^j$, and put
\begin{align*}
c_{i+1} &= c_i + r_0 \\
P_{i+1} &= P_i + \sum_{j>0} r_j T^j \\
Q_{i+1} &= Q_i + \sum_{j<0} r_j T^j.
\end{align*}
Suppose that $\min\{v(c-1), v_r(P_i-1), v_r(Q_i - 1)\} \geq v_r(R-1)$.
By Lemma~\ref{L:scale product}, $v_r(R- c_i P_i Q_i) \geq v_r(R-1)$, and
\[
v_r(R - c_{i+1} P_{i+1} Q_{i+1}) \geq v_r(R - c_i P_i Q_i) + v_r(R-1).
\]
It follows that $c_i, P_i, Q_i$ converge to 
limits $c,P,Q$ with the desired properties.
\end{proof}

\begin{cor} \label{C:one slope}
If $R(T) \in F \{T^{\pm} \}$ is irreducible, then it has only one slope.
\end{cor}

\subsection{Classification of pure $\phi$-modules}
\label{subsec:classif}

We next classify the $\phi$-modules over $k((u^\QQ))$, then
classify the pure $\phi$-modules over $\tcalR$ (Proposition~\ref{P:pure
equiv2}).

\begin{notation}
Throughout this subsection only, write $F = k((u^\QQ))$;
note that this is consistent with Notation~\ref{N:valued field} if we
put $s = q$,
take $v_F$ to be the $u$-adic valuation, and take $\phi_F$ of the form
$\sum c_i u^i \mapsto \sum \phi_k(c_i) u^{qi}$ for some automorphism
$\phi_k$ of $k$.
\end{notation}

\begin{lemma} \label{L:slope 0 root}
Let $P(T) \in F\{T\}$ be a twisted polynomial over $F$ 
with all Newton slopes equal to $0$. Then there exists
$x \in \gotho_F^*$ such that $P(\phi)(x) = 0$.
\end{lemma}
\begin{proof}
We may assume that $P$ has constant term 1.
Since $\phi$-modules over $k$
are trivial (by Hypothesis~\ref{H:field size}), we can find
$z \in \gotho_F^*$ with $P(\phi)(z) \in \gothm_F$.
Since $(P-1)(\phi)$ is contractive towards 0 on $\gothm_F$, we
can find $y \in \gothm_F$ such that $P(\phi)(y) = P(\phi)(z)$.
Put $x = z - y$; then $P(\phi)(x) = 0$.
\end{proof}

\begin{lemma} \label{L:factor slope 0}
Let $P(T) \in F\{T\}$ be a monic twisted polynomial over $F$ 
with all Newton slopes equal to $0$. Then $P(T)$ factors as a 
product $\prod_j (T - a_j)$ for some $a_j \in \gotho_F^*$.
\end{lemma}
\begin{proof}
By Lemma~\ref{L:slope 0 root}, there exists $x \in \gotho_F^*$
such that $P(\phi)(x) = 0$. By the division algorithm for
twisted polynomials,
$P(T)$ is right divisible by $T - a$ for $a = \phi(x)/x$; the claim
then follows by induction.
\end{proof}

\begin{lemma} \label{L:irreducible over k}
Every irreducible $\phi$-module over $F$ 
is trivial.
\end{lemma}
\begin{proof}
Let $V$ be an irreducible $\phi$-module over $F$;
we can then write $V$ as $F\{T^{\pm}\}/F\{T^{\pm}\}P$ for some
monic irreducible twisted polynomial $P(T)$.
By Corollary~\ref{C:one slope},
$P$ has only one slope, which we can force to be 0 by rescaling. 
By Lemma~\ref{L:factor slope 0}, $P$ must equal
$T - a$ for some $a \in \gotho_F^*$. But the equation $\phi(x) = ax$ has 
a solution $x \in \gotho_F^*$ by Lemma~\ref{L:slope 0 root}, yielding
the triviality of $V$.
\end{proof}

\begin{prop} \label{P:essential over k}
Every $\phi$-module over $F = k((u^\QQ))$ is trivial.
\end{prop}
\begin{proof}
Any $\phi$-module over $F$ can be written as a
successive extension of irreducibles, which are all trivial by
Lemma~\ref{L:irreducible over k}. 
By Corollary~\ref{C:kernel cokernel}, the extensions between trivial
$\phi$-modules all split, yielding the claim.
\end{proof}

\begin{defn}
For $P(T) = \sum_i a_i T^i 
\in F\{T^{\pm}\}$ nonzero and $z \in F$, define the \emph{inhomogeneous
Newton polygon} of the pair $(P,z)$ as the lower convex hull of the set
\[
\{(-q^i, v_F(a_i)): i \in \ZZ\} \cup \{(0, v_F(z))\};
\]
note that any slope of this polygon not involving the point $(0, v_F(z))$ is
equal to $q-1$ times a slope of the homogeneous Newton polygon.
\end{defn}

\begin{prop} \label{P:inhomogeneous}
Given $P(T) \in F\{T^{\pm}\}$ nonzero and $z\in F$, for each $r \in \RR$ occurring
as a slope of the inhomogeneous Newton polygon of $(P,z)$, there exists
$x \in F$ with $v_F(x) = r$ such that $P(\phi)(x) = z$.
\end{prop}
\begin{proof}
By applying Proposition~\ref{P:factor twisted}, we may reduce to the case
where $P$ has a single homogeneous Newton slope; by twisting, we may force 
that slope to be 0. By Lemma~\ref{L:factor slope 0}, we may reduce to
the case $P(T) = T - a$ for $a \in \gotho_F^*$. By Lemma~\ref{L:slope 0 root},
we may assume that $a = 1$; in this case, the claim follows from
Corollary~\ref{C:kernel cokernel}.
\end{proof}

Before proving Proposition~\ref{P:pure equiv2}, 
we need one
more calculation, which includes Proposition~\ref{P:cherbonnier1} 
(see also Remark~\ref{R:cherbonnier1}).

\begin{prop} \label{P:cherbonnier2}
Let $\tilde{\calE}$ denote the $\gothm_K$-adic completion of
$\tcalR^{\bd}$. Let $A$ be an $n \times n$ matrix over $\tcalR^{\inte}$. 
If $\bv
\in \tilde{\calE}^n$ is a column vector such that
$A \bv = \phi (\bv)$, then $\bv \in (\tcalR^{\bd})^n$.
\end{prop}
\begin{proof}
By rescaling by a factor of $u$ (as in the proof of
Proposition~\ref{P:h1 both ways}), we may reduce to the case where
the entries of $A$ are bounded by 1 under $|\cdot|_r$;
we may also assume $\bv$ has entries in the completion of 
$\tcalR^{\inte}$.
Write $\bv = \sum_{j=1}^n \sum_{i \in \QQ} c_{ij} u^i \be_j$,
where $\be_1, \dots, \be_n$ are the standard basis vectors;
it suffices to show that $|c_{ij}u^i|_r \leq 1$ for all $i,j$,
as then $\bv$ will have entries in $\tcalR^s$ for any $s \in (0,r)$.

Suppose the contrary; note that 
$|c_{ij}| \leq 1$ for all $i,j$ by our normalization of $\bv$,
so any pair $i,j$ with $|c_{ij} u^i|_r > 1$ must have $i<0$,
and hence 
\begin{equation} \label{eq:cherbonnier2}
|\phi^{-1}(c_{ij} u^i)|_r = |c_{ij} u^{i/q}|_r < |c_{ij} u^{i}|_r.
\end{equation}
Let $h$ be the maximum of $|c_{ij}|$ over all pairs
$i,j$ with $|c_{ij} u^i|_r > 1$. Then there is a pair $(i_0,j_0)$
with $|c_{i_0,j_0}| = h$ which maximizes $|c_{i_0,j_0} u^{i_0}|_r$.
However, if we expand $A \bv = 
\sum_{j=1}^n \sum_{i \in \QQ} d_{ij} u^i \be_j$, 
then for each pair $i,j$ with $|d_{ij}| = h$,
we have $|\phi^{-1}(d_{ij} u^i)|_r < |c_{i_0,j_0} u^{i_0}|_r$
by \eqref{eq:cherbonnier2}. This contradicts the equality
$\bv = \phi^{-1}(A \bv)$, proving the claim.
\end{proof}

We now prove Proposition~\ref{P:pure equiv2}: the categories of
pure $\phi$-modules over $K$ and over $\tcalR$ of a given slope $s$
are equivalent.
\begin{proof}[Proof of Proposition~\ref{P:pure equiv2}]
We first check full faithfulness.
By Lemma~\ref{L:pure push} and twisting, it suffices to check this for $s=0$;
that is, we must check that given an
\'etale $\phi$-module $M_0$ over $K$, we must have
$H^0(M_0) \cong H^0(M_0 \otimes_K \tcalR)$.
By Hypothesis~\ref{H:field size}, we may assume that $M_0$ is trivial;
then Lemma~\ref{L:where fixed} yields the claim.

We next check essential surjectivity; 
we may proceed as in the proof of Theorem~\ref{T:pure equiv}
to reduce to the case $s=0$. Let $M$ be an \'etale $\phi$-module
over $\tcalR$, and choose an \'etale lattice $M_0$ of $M$.
By repeated application of 
Proposition~\ref{P:essential over k}, after tensoring with the
$\gothm_K$-adic completion of $\tcalR^{\inte}$, we can find a basis of $M_0$
fixed by $\phi$. By Proposition~\ref{P:cherbonnier2}, this basis is in fact
contained in $M_0$ itself, yielding the claim.
\end{proof}

\subsection{The local calculation}
\label{subsec:local calc}

We now perform
the explicit calculation that proves Proposition~\ref{P:local
calc}, thus completing the proof of Theorem~\ref{T:semistable pure}. To avoid
notational overload, we elide a few routine calculations that can be found
in \cite{kedlaya-local}. (Also compare \cite[\S 9,10]{hartl-pink}.)

\begin{defn}
Let $\tcalR^{\tr}$ (for ``truncated'') denote the set of elements of $\tcalR$
whose support is bounded below. This forms a subring of $\tcalR$ carrying
a $u$-adic
valuation $v$. Note that a unit in $\tcalR^{\tr}$ is precisely an element
$x = \sum_i a_i u^i$ for which the support of $x$ has a least element $j$,
and for which $|a_i| \leq |a_j|$ for all $i \in \QQ$; in particular,
such elements belong to $\tcalR^{\bd}$, so we can apply 
the valuation $w$ to them.
\end{defn}

\begin{lemma} \label{L:positioning}
Let $P$ be a $\phi$-module over $K$ of rank $1$ and degree $n > 0$, and
fix a generator $\bv$ of $P$. 
\begin{enumerate}
\item[(a)]
For any $x \in \tcalR^{\tr}$ with support in $[0, +\infty)$,
the class of $x \bv$ in $H^1(P \otimes \tcalR)$ vanishes.
\item[(b)]
Each class in $H^1(P \otimes \tcalR)$ has a representative of the form
$\sum_{j=0}^{n-1} u_j \bv$, where for each $j$, either $u_j = 0$,
or $u_j \in (\tcalR^{\tr})^*$, $w(u_j) = j$, and $v(u_j) < 0$.
\end{enumerate}
\end{lemma}
\begin{proof}
For (a), we first use Hypothesis~\ref{H:field size} to eliminate
constant terms, then note that if $x$ has no constant term,
the sum $\sum_{i=0}^\infty \phi^i(x \bv)$ converges
and its limit $\bw$ satisfies $\bw - \phi(\bw) = x \bv$.
We deduce (b) from (a) plus a direct calculation;
see also \cite[Lemmas~4.13 and~4.14]{kedlaya-local}
or \cite[Lemma~4.3.2]{kedlaya-slope}.
\end{proof}

We now prove Proposition~\ref{P:local calc}: if 
$N'$ is a pure $\phi^n$-module over $\tcalR$ of rank
1 and degree 1, $P$ is a pure $\phi$-module over $\tcalR$ of rank 1 and
degree -1, and
\begin{equation} \label{eq:exact seq}
0 \to [n]^* N' \to M \to P \to 0
\end{equation}
is a short exact sequence of $\phi$-modules, then $H^0(M) \neq 0$.

\begin{proof}[Proof of Proposition~\ref{P:local calc}]
The snake lemma gives an exact sequence
\[
H^0(M) \to H^0(P) \to H^1([n]^* N'),
\]
where the second map is pairing with the class $\alpha \in H^1(P^\dual
\otimes [n]^* N')$ corresponding to the extension \eqref{eq:exact seq};
it suffices to show that this second map has nonzero kernel.

Note that $P^\dual \otimes [n]^* N' \cong [n]^* ([n]_* P^\dual \otimes N')$
as in Definition~\ref{D:pull push}, so we may view $\alpha$ as an element of
$H^1([n]^*([n]_* P^\dual \otimes N')) \cong H^1([n]_* P^\dual \otimes N')$.
Similarly,
we may view the pairing with $\alpha$ as the composition of the map
$H^0(P) \to H^0([n]_* P)$ with the map $H^0([n]_* P) \to 
H^1(N')$ given by pairing with the class in $H^1([n]_* P^\dual \otimes N')$.
If the class vanishes, there is nothing to check, so we may assume that
it does not vanish.

By Proposition~\ref{P:pure equiv2}, $P$ and $N'$ are obtained by base
change from certain $\phi$- and $\phi^n$-modules $P_0$ and $N'_0$,
respectively, over $K$;
choose generators $\bv$ and $\bw$ of $P_0$ and $N'_0$, and define
$\lambda, \mu \in K^*$ by $\phi(\bv) = \lambda \bv$ and $\phi^n(\bw)
= \mu \bw$.
Put $Q_0 = [n]_* P_0^\dual \otimes N'_0$ and $Q = [n]_* P^\dual \otimes N'
\cong Q_0 \otimes_K \tcalR$; let $\bx$ be the generator $\bv^\dual \otimes \bw$
of $Q_0$ (where $\bv^\dual$ is the generator of $P^\dual$ dual to $\bv$).
By Lemma~\ref{L:positioning}, we can then represent
the class $\alpha \in H^1(Q)$ by a nonzero element of $Q$ of the form
$\sum_{j=0}^n u_j \bx$, where each
$u_j$ is either zero or a unit in $\tcalR^{\tr}$ with $w(u_j) = j$
and $v(u_j) < 0$.

We now follow \cite[Lemma~4.12]{kedlaya-local}.
For $j \in \{0, \dots, n\}$ such that $u_j \neq 0$, 
$l \in \ZZ$, and $m \in (0, +\infty)$, define
\[
e(j,l,m) = (v(u_j) + m q^{-l})q^{-n(j+l)}.
\]
For fixed $j$ and $m$, $e(j,l,m)$ approaches 0 from below as $l \to +\infty$,
and tends to $+\infty$ as $l \to -\infty$. Hence the minimum
$h(m) = \min_{j,l} \{e(j,l,m)\}$ is well-defined; we observe that
$h$ is a continuous, piecewise linear, and 
increasing map from
$(0,+\infty)$ to $(-\infty,0)$, and that
$h(qm) = q^{-n} h(m)$ because $e(j,l+1,qm) = q^{-n} e(j,l,m)$. Another interpretation
is that the lower convex hull of the set $H$ of points
\[
(-q^{-nj-(n+1)l}, q^{-nj-nl} v(u_j)) \qquad (j=0,\dots,n; \quad l \in \ZZ)
\]
has all slopes positive, and all segments finite.

Pick $r \in (0,+\infty)$ at which $h$ changes slope; that is,
$r$ is a slope of the lower convex hull of $H$. Let $S$ denote the set
of ordered pairs $(j,l)$ for which $e(j,l,r) < q^{-n} h(r)$; this set
is finite. Let $T$ be the set of ordered pairs $(j,l)$ for which
$e(j,l,r) < 0$; this set (which contains $S$) is infinite, but the values
of $l$ for pairs $(j,l) \in T$ are bounded below. For each pair $(j,l)$,
put $s(j,l) = \lfloor \log_{q^n}(h(r)/e(j,l,r)) \rfloor$. Then
the following properties hold.
\begin{enumerate}
\item[(a)]
For $(j,l) \in T$, $s(j,l) \geq 0$.
\item[(b)]
For $(j,l) \in T$, $e(j,l,r) q^{ns(j,l)} \in [h(r), q^{-n}h(r))$.
\item[(c)]
We have $(j,l) \in S$ if and only if $(j,l) \in T$ and $s(j,l) = 0$.
\item[(d)]
For any $c > 0$, there are only finitely many pairs $(j,l) \in T$ with
$s(j,l) \leq c$.
\end{enumerate}

Define the twisted powers $\lambda^{\{m\}}$ and $\mu^{\{m\}}$ of $\lambda$
and $\mu$ by the two-way recurrences
\begin{gather*}
\lambda^{\{0\}} = 1, \qquad
\lambda^{\{m+1\}} = \phi(\lambda^{\{m\}}) \lambda \\
\mu^{\{0\}} = 1, \qquad
\mu^{\{m+1\}} = \phi^n(\mu^{\{m\}}) \mu.
\end{gather*}

For $c \in \RR$, let $U_c$ be the set of $z \in \tcalR^{\tr} \cap \tcalR^{\inte}$
with $v(z) \geq c$. Then the function
\[
R(z) = \sum_{(j,l) \in T} 
\mu^{\{-j-l+s(j,l)\}} \phi^{-nj-nl+ns(j,l)} (u_j \lambda^{\{-l\}} \phi^{-l}(z))
\]
carries $U_r$ into $U_{h(r)}$ by a direct calculation.
Modulo $\pi$, we have
\begin{equation} \label{eq:puiseux poly}
R(z) \equiv \sum_{(j,l) \in S} 
\mu^{\{-j-l\}} \phi^{-nj-nl} (u_j \lambda^{\{-l\}} \phi^{-l}(z));
\end{equation}
note that the values $-nj-(n+1)l$ are distinct for all $(j,l) \in S$,
since $j$ only runs over $\{0,\dots,n\}$.
Write the reduction modulo $\pi$
of the right side of \eqref{eq:puiseux poly} as
$Q(\phi)(z)$ for some twisted Laurent polynomial
$Q(T) \in F\{T^{\pm}\}$ with $F = k((u^{\QQ}))$.
By Proposition~\ref{P:inhomogeneous} applied repeatedly, we
can construct a nonzero $z \in U_r$ such that $R(z) = 0$.

One now calculates using Lemma~\ref{L:positioning}(a)
(see \cite[Lemma~4.12]{kedlaya-local} for the full calculation) that
the element 
\[
\sum_{l \in \ZZ} \phi^{-l}(z \bv) = \sum_{l \in \ZZ}
\lambda^{\{-l\}} \phi^{-l}(z) \bv \in H^0(P)
\]
pairs to
zero with the class of $\alpha$. This yields the desired result.
\end{proof}

\section{Descending the slope filtration}
\label{sec:descend}

As noted at the beginning of the previous section, the proof of the slope
filtration theorem (Theorem~\ref{T:slope filt})
consists of two stages, the first of which (classifying
$\phi$-modules over the overring $\tcalR$ of $\calR$) has been accomplished in 
the previous section. In this section, we explain how to descend the
resulting slope filtration from $\tcalR$ back to $\calR$.

As was done in the previous section, we recommend on a first reading 
to read only the overview (Subsection~\ref{subsec:descend overview}),
then return later for the technical details.

\subsection{Overview}
\label{subsec:descend overview}

\begin{defn} \label{D:descent setup}
We now revert to allowing $K$ to be an arbitrary field as in
Definition~\ref{D:initial}. Choose a complete extension $L$ of $K$
with the same value group, 
admitting an extension $\phi$ to an automorphism, such that
every \'etale $\phi$-module
over $L$ is trivial. More precisely, form such an $L$ by first
taking the completed direct limit of $K \stackrel{\phi}{\to} 
K \stackrel{\phi}{\to} \cdots$ and then applying Proposition~\ref{P:splitting
field} below.
Under these conditions, we can embed $\calR_K$ into $\calR_L$, and then
embed $\calR_L$ into $\tcalR_L$ as in Proposition~\ref{P:embedding}.
\end{defn}

Recall that we are trying to prove Theorem~\ref{T:slope filt}, which
states that every module-semistable $\phi$-module over $\calR$ is pure.
As noted earlier, this result 
follows from Theorem~\ref{T:semistable pure} (which asserts
that module-semistable $\phi$-modules over $\tcalR_L$ are pure)
plus the following assertions.

\begin{theorem} \label{T:ascend semistable}
Let $M$ be a module-semistable $\phi$-module over $\calR$. Then
$M \otimes \tcalR_L$ is module-semistable.
\end{theorem}

\begin{theorem} \label{T:descend pure}
Let $M$ be a $\phi$-module over $\calR$ such that
$M \otimes \tcalR_L$ is pure. Then $M$ is pure.
\end{theorem}

The proofs of Theorems~\ref{T:ascend semistable} and~\ref{T:descend pure}
 amount to faithfully flat descent:
Theorem~\ref{T:ascend semistable} relies on the fact that the first
step of the HN filtration of $M \otimes \tcalR_L$ descends to $\calR$,
while
Theorem~\ref{T:descend pure} depends on the fact that the
pure $\phi$-module over $\tcalR_L^{\bd}$ obtained by descending
$M \otimes \tcalR_L$ itself descends to $\calR^{\bd}$.
The rest of this section will be occupied with setting up the descent
formalism, then making the calculations that allow the use of faithfully
flat descent.

\subsection{Splitting \'etale $\phi$-modules}

We now construct the field $L$ demanded by Definition~\ref{D:descent setup}.

\begin{defn}
Suppose that $\phi_K$ is bijective.
By an \emph{admissible extension} of $K$, we will mean a field $L$ containing 
$K$, complete for a nonarchimedean absolute value extending the one on $K$
with the same value group, and equipped with an isometric field automorphism
$\phi_L$ extending $\phi_K$.
\end{defn}

\begin{lemma} \label{L:split h1}
For any $z \in K^*$,
there exists an admissible extension $L$ of $K$ such that the equation
$\phi(x) - x = z$ has a solution $x \in L$.
\end{lemma}
\begin{proof}
Let $L$ be the completion of the rational function field $K(x)$ for the
Gauss norm with $|x| = |z|$. Extend $\phi_K$ to an automorphism $\phi_L$ of
$L$ by setting
$\phi_L(x) = x + z$.
\end{proof}

\begin{lemma} \label{L:split poly}
Let $P(T) = T^n + a_{n-1} T^{n-1} + \cdots + a_0$ be a twisted polynomial
over $\gotho_K$ with $|a_0| = 1$. Then there exists an admissible extension
$L$ of $K$ such that the equation $P(\phi)(x) = 0$ has a
solution $x \in \gotho_L^*$.
\end{lemma}
\begin{proof}
Let $L$ be the completion of the rational function field 
$K(y_0, \dots, y_{n-1})$
under the Gauss norm normalized with $|y_0| = \cdots = |y_{n-1}| = 1$.
Extend $\phi_K$ to an automorphism $\phi_L$ of $L$ by setting
$\phi_L(y_i) = y_{i+1}$ for $i=0, \dots, n-2$ and
$\phi_L(y_{n-1}) = -a_{n-1}y_{n-1} - \cdots - a_0 y_0$, then
take $x = y_0$.
\end{proof}

\begin{prop} \label{P:splitting field}
There exists a complete extension $L$ of $K$ with the same value group,
equipped with an extension of $\phi_K$, such that 
any \'etale $\phi$-module over $L$ is trivial.
\end{prop}
\begin{proof}
It suffices to construct $L$ 
trivializing a single irreducible \'etale $\phi$-module
$M$ over $K$, as we can construct the desired field by transfinitely iterating
this construction and completing at all limit stages. 

Since $M$ is irreducible, we must have $M \cong K\{T^{\pm}\}/K\{T^{\pm}\}P(T)$
for some irreducible monic twisted polynomial $P(T)$. If we write
$P(T) = T^n + a_{n-1} T^{n-1} + \cdots + a_0$, then $|a_0| = 1$ because $\deg(M) = 0$.
By Corollary~\ref{C:one slope} (in the case $s=1$), 
$P$ can only have one Newton slope, which must be 0;
hence $P(T)$ has coefficients in $\gotho_K$. We can then apply
Lemma~\ref{L:split poly} to construct $L$ over which the equation
$P(\phi)(x) = 0$ has a solution $x \in \gotho_L^*$; 
that solution
gives rise to a nontrivial $\phi$-submodule of $M$. 

Repeating
the construction, we obtain a field over which $M$ becomes a successive
extension of trivial \'etale $\phi$-modules of rank 1. By repeated use of
Lemma~\ref{L:split h1}, we can split this filtration 
by passing to a suitably large $L$.
This yields the claim.
\end{proof}

\begin{remark}
Note that the field $L$ constructed above is not a Picard-Vessiot extension
of $K$ in the sense of the Galois theory of difference fields;
this Galois theory is a bit
complicated because it cannot be carried out within the category of fields,
as examples like the difference equation $\phi(x) = -x$ show. 
See \cite[Chapter~1]{singer-vanderput} for more discussion of this point,
and a development of difference Galois theory in a restricted setting;
see also \cite{andre-diff} for a more general development. (Thanks
to Michael Singer for pointing out this reference.)
\end{remark}

\subsection{The use of faithfully flat descent}

In this subsection, we set up faithfully flat descent and illustrate
how we will use it to prove
Theorems~\ref{T:ascend semistable}
and~\ref{T:descend pure}.

\begin{defn} \label{D:descent}
Let $R \to S$ be a faithfully flat morphism of rings
equipped with
compatible endomorphisms $\phi$. Let $M$ be a $\phi$-module over $R$,
put $M_S = M \otimes_R S$, and let $N_S$ be a $\phi$-submodule of $M_S$.
We say that $N_S$ \emph{descends to $R$} if there exists a
$\phi$-submodule $N$ of $M$ such that the image of
$N \otimes_R S$ in $M_S$ coincides with $N_S$.
We say a filtration descends to $R$ if each 
term does so.
\end{defn}

\begin{prop} \label{P:descent}
Let $R \to S$ be a faithfully flat morphism of domains
equipped with compatible endomorphisms $\phi$.
Put $S_2 = S \otimes_R S$ and define $i_1, i_2: S \to S_2$ by
$i_1(s) = s \otimes 1$ and $i_2(s) = 1 \otimes s$.
Let $M$ be a $\phi$-module over $R$, put $M_S = M \otimes_R S$, and let
$N_S$ be a $\phi$-submodule of $M_S$.
Then $N_S$ descends to $R$ if and only if $N \otimes_{i_1} S_2 =
N \otimes_{i_2} S_2$ within $M \otimes_R S_2$; moreover, if this occurs,
then there is a \emph{unique} $\phi$-submodule $N$ of $M$ such that
$N_S = N \otimes_R S$ within $M_S$.
\end{prop}
\begin{proof}
The equality $N \otimes_{i_1} S_2 =
N \otimes_{i_2} S_2$ implies that the effective descent datum obtained from $M$
induces a descent datum on $N$ (the cocycle condition can be checked on $M$).
We may thus apply faithfully flat descent for modules
\cite[Expos\'e~VIII, Corollaire~1.3]{grothendieck} to conclude.
\end{proof}

We use faithfully flat descent as follows.
\begin{defn}
Define 
\begin{align*}
\calS &= \tcalR_L \otimes_{\calR} \tcalR_L \\
\calS^{\bd} &= \tcalR_L^{\bd} \otimes_{\calR^{\bd}} \tcalR_L^{\bd} \\
\calS^{\inte} &= \tcalR_L^{\inte} \otimes_{\calR^{\inte}} \tcalR_L^{\inte}.
\end{align*}
We will show later that $\calR \to \tcalR_L$, $\calR^{\bd} \to \tcalR_L^{\bd}$
are faithfully flat and that
$\calS^{\bd} \to \calS$ is injective
(Remark~\ref{R:faithfully flat}).
\end{defn}

The following weak analogue of Proposition~\ref{P:h1 both ways}
will be proved in Subsection~\ref{subsec:tensor}.
\begin{prop} \label{P:tensor bounded}
Let $A$ be an $n \times n$ matrix over $\calS^{\inte}$, and let
$\bv$ be a column vector over $\calS$ 
such that $\bv = A \phi(\bv)$. Then
$\bv$ has entries in $\calS^{\bd}$.
\end{prop}

We now demonstrate how Proposition~\ref{P:tensor bounded} can be used
to establish the theorems asserted at the start of this section.

\begin{proof}[Proof of Theorem~\ref{T:ascend semistable}]
Suppose that $M \otimes \tcalR_L$ is not semistable.
Let $0 = M_{L,0} \subset M_{L,1} \subset \cdots \subset M_{L,l} = M_L$
denote the HN filtration of $M_L = M \otimes \tcalR_L$. We will show
that $M_{L,1} \otimes_{i_2} S_2 \subseteq M_{L,j} \otimes_{i_1} S_2$ for
$j = l,l-1,\dots,1$ by descending induction; the base case $j=l$ is trivial.

Given that $M_{L,1} \otimes_{i_2} S_2 \subset M_{L,j} \otimes_{i_1} S_2$
for some $j > 1$,
we get a homomorphism 
\[
M_{L,1} \otimes_{i_2} S_2 \to (M_{L,j}/M_{L,j-1}) \otimes_{i_1} S_2.
\]
Since $M_{L,1}$ and $M_{L,j}/M_{L,j-1}$ are pure and
$\mu(M_{L,1}) < \mu(M_{L,j}/M_{L,j-1})$, this homomorphism is forced to
vanish: otherwise, 
by Proposition~\ref{P:tensor bounded} the morphism would be defined
over $\calS^{\bd}$, but in that case it would have to preserve slopes
because $\calS^{\bd}$ carries an $\gothm_K$-adic valuation.
Hence $M_{L,1} \otimes_{i_2} S_2 \subseteq M_{L,j-1} \otimes_{i_1} S_2$,
completing the induction.

The induction shows that $M_{L,1}$ satisfies the condition for
faithfully flat descent (Proposition~\ref{P:descent}), so it descends
to $\calR$. Hence $M$ cannot be semistable either.
\end{proof}

\begin{proof} [Proof of Theorem~\ref{T:descend pure}]
By applying $[a]_*$ (invoking 
Lemma~\ref{L:pure push})
and twisting, we may reduce to the case $\mu(M) = 0$,
so $M \otimes \tcalR_L$ is \'etale.
Choose a basis $\bv_1, \dots, \bv_n$ of an \'etale lattice of 
$M \otimes \tcalR_L$, so that the matrix
$A$ defined by $\phi(\bv_j) = \sum_i A_{ij} \bv_i$ is 
invertible over $\tcalR_L^{\inte}$.

There exists an invertible change-of-basis matrix $U$
over $\calS$ such that
\[
\bv_j \otimes_{i_1} 1 = \sum_i U_{ij} (\bv_i \otimes_{i_2} 1).
\]
Upon applying $\phi$ to both sides, we deduce that
$U(A \otimes_{i_1} 1) = (A \otimes_{i_2} 1)\phi(U)$.
By Proposition~\ref{P:tensor bounded},
$U$ has entries in $\calS^{\bd}$, as does its inverse by the same
argument with $M$ replaced by $M^\dual$. 
Hence by Proposition~\ref{P:descent}, $M$ descends to $\calR^{\bd}$;
let $N$ be the resulting $\phi$-module over $\calR^{\bd}$.

Choose any basis of $N$ and let $P$ be the $\calR^{\inte}$-span of the images
of the basis elements under powers of $\phi$. 
By computing in terms of $\bv_1, \dots,
\bv_n$, we see that $P$ is bounded, hence is a $\phi$-stable 
$\calR^{\inte}$-lattice
in $M$. By Lemma~\ref{L:etale lattice}, $P \otimes \tcalR_L^{\inte}$
is a $\phi$-module, as then must be $P$. Thus $M$ is \'etale, as desired.
\end{proof}

It now remains to prove the faithful flatness results and to make the
calculation to check Proposition~\ref{P:tensor bounded}; these occupy the
remainder of the chapter.

\subsection{Interlude: tensoring over B\'ezout domains}

In order to use faithfully flat descent for our purposes, 
it will help to gather a few facts about
tensoring over B\'ezout domains.

\begin{prop} \label{P:faithfully flat1}
Let $R \hookrightarrow S$ be an inclusion of domains with $R$ 
B\'ezout. Then $S$ is faithfully flat over $R$ if and only if
$S^* \cap R = R^*$.
\end{prop}
\begin{proof}
Recall that $S$ is flat (resp.\ faithful) over $R$ if and only if for each
finitely generated proper ideal $I$ of $R$, the multiplication
map $I \otimes S \to S$ is injective (resp.\ not surjective).
Since $R$ is B\'ezout, $I$ admits a single generator $r \notin R^*$, and
$I \otimes S = rR \otimes S \cong rS$, so the map $I \otimes S \to S$
is injective, and it is surjective if and only if $r \in S^*$.
This yields the claim.
\end{proof}

\begin{lemma} \label{L:rewrite tensor}
Let $M,N$ be modules over a B\'ezout domain $R$. Given a presentation
$\sum_{i=1}^n y_i \otimes z_i$ of $x \in M \otimes_{R} N$ and elements
$u_1, \dots, u_n \in R$ generating the unit ideal,
there exists another presentation
$\sum_{j=1}^n y'_j \otimes z'_j$ of $x$ with $y'_1 = \sum_{i=1}^n u_i y_i$.
\end{lemma}
\begin{proof}
By  \cite[Lemma~2.3]{kedlaya-local},
we can construct an invertible matrix $U$ over $R$ with $U_{i1} = u_i$
for $i=1,\dots,n$.
Then
\begin{align*}
\sum_i y_i \otimes z_i &=
 \sum_{i,j,l} U_{ij} (U^{-1})_{jl} y_i \otimes z_l \\
&= \sum_j \left( \sum_i U_{ij} y_i \right) \otimes
\left( \sum_l (U^{-1})_{jl} z_l \right),
\end{align*}
so we may take $y'_j = \sum_{i=1}^n U_{ij} y_i$ and $z'_j = \sum_{l=1}^n (U^{-1})_{jl} z_l$.
\end{proof}

\begin{cor} \label{C:lin ind}
Let $M,N$ be modules over a B\'ezout domain $R$. 
If $\sum_{i=1}^n y_i \otimes z_i$ is a presentation of 
some $x \in M \otimes_R N$ with $n$ minimal,
then $y_1, \dots, y_n$ are linearly independent over $R$.
\end{cor}
\begin{proof}
If on the contrary $y_1, \dots, y_n$ are linearly dependent over $R$,
then we can find $u_1, \dots, u_n \in R$ such that $u_1 y_1 + \cdots + u_n y_n = 0$.
By the B\'ezout property, $u_1, \dots, u_n$ generate a principal ideal, so we
can divide through by a generator to reduce to the case where $u_1, \dots, u_n$
generate the unit ideal. Applying Lemma~\ref{L:rewrite tensor} now yields
a contradiction to the minimality of $n$.
\end{proof}

\subsection{Projections}
\label{subsec:tensor}

The key to the descent argument is the construction of a certain projection
from $\tcalR_L$ back to $\calR$, sectioning the inclusion going the other way
that was constructed by Proposition~\ref{P:embedding}.
We now construct this projection, then use it to resolve all the outstanding
statements needed to complete the proof of Theorem~\ref{T:slope filt}.

\begin{defn}
Let $\ell$ be the residue field of $L$,
fix a basis $\overline{B}$ of $\ell$ over $k$ containing 1,
lift $\overline{B}$ to a subset $B$ of $\gotho_L$ containing 1,
and fix a uniformizer $\pi$ of $K$.
Then as in \cite[Proposition~4.1]{kedlaya-full}, one sees that every
element $x \in \tcalR^{\inte}_L/\gothm_K^n \tcalR_L^{\inte}$ can be written
uniquely as a formal sum 
\[
\sum_{\alpha \in [0,1) \cap \QQ} \sum_{b \in B} x_{\alpha,b} u^\alpha b
\qquad (x_{\alpha,b} \in \calR^{\inte}/\gothm_K^n \calR^{\inte})
\]
in which:
\begin{itemize}
\item
for each $\alpha 
\in [0,1) \cap \QQ$, there are only finitely many $b$ for which $x_{\alpha,b} 
\neq 0$;
\item
if we write $S_c$ for the
set of $\alpha \in [0,1) \cap \QQ$ for 
which the $t$-adic valuation of any $x_{\alpha,b}$ 
(which is well-defined because
$x_{\alpha,b}$ is truncated modulo $\pi^n$) is less than $c$, then
$S_c$ is well-ordered for all $c$ and empty for sufficiently small $c$.
\end{itemize}
Given $x$ thusly presented, write $f(x) = x_{0,1}$; then
again as in \cite[Proposition~4.1]{kedlaya-full},
one checks that  for $r_0$ as in Remark~\ref{R:annuli} and $r \in (0,r_0)$,
$f$ induces a continuous map $\tcalR_L^r \to \calR^r$ with the
property that for $x \in \tcalR_L^r$,
\begin{equation} \label{eq:compare bounds}
|x|_r = \sup_{\alpha \in [0,1) \cap \QQ, a \in L^*} \{
|a|^{-1} e^{-\alpha r} |f(a u^{-\alpha} x)|_r\}.
\end{equation}
(Compare also 
\cite[Proposition~8.1]{dejong} and \cite[Lemma~2.2.19]{kedlaya-slope}.)
\end{defn}

\begin{prop} \label{P:mult}
The multiplication map $\tcalR_L^{\bd} \otimes_{\calR^{\bd}} \calR \to
\tcalR_L$ is injective.
\end{prop}
\begin{proof}
Suppose the contrary; choose $x \neq 0$ 
in the kernel of the multiplication map, and
choose a presentation $x = \sum_{i=1}^n y_i \otimes z_i$ with $n$ minimal.
Then $z_1, \dots, z_n$ are linearly independent over $\calR^{\bd}$ by
Corollary~\ref{C:lin ind}. On the other hand, as a
corollary of \eqref{eq:compare bounds}, we may choose $\alpha \in
[0,1) \cap \QQ$ and $a \in L^*$ such that
$f(a u^{-\alpha} y_1) \neq 0$; we then obtain the nontrivial dependence
relation
$0 = \sum_{i=1}^n f(a u^{-\alpha} y_i) z_i$, contradiction.
\end{proof}

\begin{remark} \label{R:faithfully flat}
We now have a number of faithfully flat inclusions.
For one, $\calR^{\bd} \to \calR$ is faithfully flat by
Proposition~\ref{P:faithfully flat1} and the fact that
$\calR^* = (\calR^{\bd})^*$ (Remark~\ref{R:same units}).
For another, $\calR \to \tcalR_L$ is faithfully flat by
Proposition~\ref{P:faithfully flat1} and the fact that
$\tcalR_L^* = (\tcalR_L^{\bd})^*$ (Remark~\ref{R:analytic ring});
similarly, $\calR^{\bd} \to \tcalR_L^{\bd}$ is faithfully flat.
Putting these together and using Proposition~\ref{P:mult}
yields injections
\[
\tcalR_L^{\bd} \otimes_{\calR^{\bd}} \tcalR_L^{\bd}
\hookrightarrow \tcalR_L^{\bd} \otimes_{\calR^{\bd}} \tcalR_L
\cong (\tcalR_L^{\bd} \otimes_{\calR^{\bd}} \calR) \otimes_{\calR} \tcalR_L
\hookrightarrow \tcalR_L \otimes_{\calR} \tcalR_L;
\]
that is, $\calS^{\bd} \to \calS$ is injective.
\end{remark}

In order to calculate on $\calS$, we use the following two-variable 
analogue of \eqref{eq:compare bounds}.
\begin{lemma} \label{L:bounded}
For $x \in \calS$, we have $x \in \calS^{\bd}$ if and only if
for some $r>0$, 
the quantities 
\begin{equation} \label{eq:bounded}
|ab|^{-1} e^{-\alpha s - \beta s} |(f \otimes f)((a u^{-\alpha} 
\otimes b u^{-\beta})x)|_s
\end{equation}
are bounded over all $s \in (0,r]$, all $a,b \in L^*$, and all
$\alpha,\beta \in [0,1) \cap \QQ$.
\end{lemma}
\begin{proof}
If $x \in \calS^{\bd}$, then we can bound the quantity
\eqref{eq:bounded} by bounding each term in a presentation of $x$.
Conversely, suppose the quantity \eqref{eq:bounded} is bounded.
Choose a presentation $x = \sum_{i=1}^n y_i \otimes z_i$
with $y_i, z_i \in \tcalR_L$ and $n$ minimal. 
We proceed by induction on $n$; we may assume
$x \neq 0$. Then $y_1 \neq 0$, so we can choose
$a, \alpha$ with $f(a u^{-\alpha} y_1) \neq 0$.

By \eqref{eq:compare bounds}, $\sum_{i=1}^n f(a u^{-\alpha} y_i) z_i 
\in \tcalR_L^{\bd}$; in particular, the ideal generated by 
the $f(a u^{-\alpha} y_i)$ in $\calR$ 
extends to the unit ideal in $\tcalR_L$. Since the ideal in $\calR$
is finitely generated, it is principal, and since
$\tcalR_L^* = (\tcalR_L^{\bd})^*$, the generator in $\calR$ must already
be a unit. That is, the $f(a u^{-\alpha} y_i)$ generate the unit ideal
in $\calR$; by Lemma~\ref{L:rewrite tensor}, we can choose another presentation
$x = \sum_{i=1}^n y'_i \otimes z'_i$ with $z'_1 = \sum_{i=1}^n 
f(a u^{-\alpha} y_i) z_i \in \tcalR_L^{\bd}$. We must have $z'_1 \neq 0$
to avoid contradicting the minimality of $n$.

Pick $b, \beta$ so that $f(b u^{-\beta} z'_1)$ is nonzero and hence is a unit
in $\calR$ (since it must lie in $\calR^{\bd}$). 
Put $c_i = f(b u^{-\beta} z'_i) / f(b u^{-\beta} z'_1)$ for $i=2, \dots, n$,
then set
\[
y''_i = 
\begin{cases}
y'_1 + c_2 y'_2 + \cdots + c_n y'_n & i = 1 \\
y'_i & i > 1,
\end{cases}
\qquad
z''_i = 
\begin{cases}
z'_i & i = 1 \\
z'_i - c_i z'_1 & i > 1,
\end{cases}
\]
so that $x = \sum_{i=1}^n y_i'' \otimes z_i''$.
Then $f(b u^{-\beta} z''_i) = 0$ for $i=2, \dots, n$,
so $y''_1 f(b u^{-\beta} z''_1) = \sum_{i=1}^n y''_i f(b u^{-\beta} z''_i) \in
\tcalR_L^{\bd}$ by \eqref{eq:compare bounds}.
Since already $f(b u^{-\beta} 
z''_1) \in
\calR^{\bd}$, we have $y''_1 \in \tcalR_L^{\bd}$. Applying the induction
hypothesis to $x - y''_1 \otimes z''_1 = \sum_{i=2}^n y''_i \otimes z''_i$ yields the claim.
\end{proof}

\begin{proof}[Proof of Proposition~\ref{P:tensor bounded}]
For each entry $\bv_i$ of $\bv$, choose a presentation
$\sum_{j} y_{ij} \otimes z_{ij}$ with $y_{ij}, z_{ij} \in \tcalR_L$.
As in the proof of Proposition~\ref{P:h1 both ways},
after possibly rescaling by a power of $u$,
we may choose $r \in (0,r_0)$ such that each term in a presentation of
$A$ has entries in $\tcalR_L^r$ and
is bounded by 1 on the annulus $e^{-r} \leq |u| < 1$;
we may also ensure that 
$y_{ij},z_{ij} \in \tcalR_L^r$ for all $i,j$.
Choose $c>0$ such that for $s \in [r/q,r]$ and all $i,j$,
$|y_{ij}|_s \leq c$ and $|z_{ij}|_s \leq c$ (possible because we are
picking $s$ in a closed interval);
then for all nonnegative integers $m$, we have
$|\phi^m(y_{ij})|_{s/q^m} \leq c$ and $|\phi^m(z_{ij})|_{s/q^m} \leq c$.
{}From the equation
\[
\bv = A \phi(A) \cdots \phi^{m-1}(A) \phi^m(\bv),
\]
we deduce that for all $\alpha,\beta \in [0,1)$ and all $a,b \in L^*$,
\[
|ab|^{-1} e^{-\alpha s-\beta s} |(f \otimes f)((a u^{-\alpha} \otimes 
b u^{-\beta}) \bv)|_s \leq c
\]
for all $s \in [r/q^{m+1}, r/q^m]$; by varying $m$, we get the same
conclusion for all $s \in (0,r]$. By Lemma~\ref{L:bounded},
$\bv$ has entries in $\calS^{\bd}$, as desired.
\end{proof}

\end{document}